\documentclass[12pt]{amsart}
\usepackage[margin=1.2in]{geometry}
\usepackage{booktabs, enumitem}
\usepackage{hyperref}
\usepackage{algorithm, algorithmicx}
\usepackage[noend]{algpseudocode}
\usepackage{tikz, subcaption, tikz-cd}
\usetikzlibrary{arrows.meta, arrows}
\usepackage[justification = centering]{caption} 
\usepackage{amsmath, mathtools}
\usepackage[bibstyle=ieee, citestyle=numeric, backend=bibtex, sortcites=false, sorting=nyt, maxnames=5]{biblatex}
\DeclareNameAlias{default}{family-given}

\bibliography{../../LateX/biblio}

\title[Projection onto the core]{Projection onto the core: An optimal reallocation to correct market failure}
\author{Dylan Laplace Mermoud}
\address{Unit{\'e} de Math{\'e}matiques Appliqu{\'e}es, ENSTA Paris, Institut Polytechnique de Paris, 91120 Palaiseau, France}
\email{dylan.laplace.mermoud@gmail.com}
\date{\today}

\newtheorem{theorem}{Theorem}

\newtheorem{proposition}{Proposition}
\newtheorem{lemma}{Lemma}

\theoremstyle{definition}
\newtheorem{definition}{Definition}
\newtheorem{example}{Example}
\newtheorem{remark}{Remark}

\AtBeginEnvironment{example}{%
  \pushQED{\qed}%
}
\AtEndEnvironment{example}{\popQED\endexample}

\DeclareMathOperator*{\argmin}{argmin}
\DeclareMathOperator*{\argmax}{argmax}

\begin{document}

\begin{abstract}
This paper provides formulae and algorithms to compute the projection onto the core of a preimputation outside it. The core of a game is described using an exponential number of linear constraints, and we cannot know beforehand which are redundant or defining the polytope. We apply these new results to market games, a class of games in which every game has a nonempty core. Given an initial state of the game represented by a preimputation, it is not guaranteed that the state of the game evolves toward the core following the dynamics induced by the domination relations. Our results identify and compute the most efficient side payment that acts on a given state of the game and yields its closest core allocation. Using this side payment, we propose a way to evaluate the failure of a market to reach a state of the economy belonging to the core, and we propose a new solution concept consisting of preimputations that minimizes this failure. 
\end{abstract}

\maketitle


\section{Introduction}

The first solution concept introduced for solving cooperative games are the \emph{stable sets}, which, according to \textcite{von1944theory}, model the standards of behavior of the social context represented by the game under consideration. The point was to ``embrace the phenomena of social equilibrium and changes'' and ``the idea of capturing social configurations'' (\textcite{leonard2010neumann}). 

\medskip 

Rapidly after, the complementary solution concept of the \emph{core} is introduced by Shapley (\textcite{zhao2018three, shubik1992game}), and its relation with the stable sets studied by \textcite{gillies1959solutions}. The core is defined as the polytope, possibly empty, of all allocations of the grand coalition revenue to all players such that no coalition can provide a better allocation for all its members. 

\medskip 

The property of having a nonempty core is an important one for a game, because ``a game that has a core has less potential for social conflict than one without a core'' (\textcite{shubik1982game}). Indeed, the core can be thought of as the set of socially fair allocations, while the stable sets represent the natural social order of the social context modeled by the game under consideration. Even for a game with a nonempty core, if this core is not a stable set, nothing ensures that the distribution of the players' endowments evolves towards a fair state. 

\medskip 

To solve this issue, we propose formulae and algorithms to compute the most efficient reallocation of players' endowments. Given any initial state represented by a preimputation, that is, an allocation that distributes the whole revenue of the grand coalition among all the players, we are able to compute the closest core allocation. This redistribution of the players' endowments does not rely on the properties of the stable sets, but on an external force or arbitrator. 

\medskip 

More formally, we compute for any preimputation outside the core its projection onto the core. According to the Hilbert Projection Theorem, the projection of an allocation outside the core onto it is the closest preimputation in the core from the initial allocation. The task of projecting onto the core is difficult, because the H-description of the core requires an exponential number of values, and we don't know in advance which of these values are relevant.

\medskip 

An important tool we use is a projector onto some specific affine subspaces of the space of allocations, for which exact formulae are provided. Moreover, we give an algorithm which associates, to any preimputation, a collection of coalitions for which the projection that gives a payment exactly equal to their value belongs to the core.   









\medskip 

Finally, we use these formulas and algorithms to compute the distance between a preimputation and the core of a game, and propose a map  which associates to each preimputation a value measuring the failure of this given preimputation to be in the core, which we apply to market games. 

\section{Preliminaries}\label{section: preliminaries}

Let \(N\) be finite set of \(n\) players. We want to know whether the players in \(N\) would agree to fully cooperate. For that to be the case, forming the coalition \(N\) must be in the interest of every player. In particular, there should not be any subcoalition of \(N\) able to provide to its players more than what \(N\) can provide to them. 

\medskip 

To formally quantify this idea, we use a \emph{coalitional function}, a set function \(v\) from the power set \(2^N\) of \(N\) to the real numbers \( \mathbb{R} \), which associates to each subset \(S\) of \(N\) a real number \(v(S)\), called the \emph{worth} of \(S\). We can interpret the worth of \(S\) as the amount of utility that the players in \(S\) have by cooperating during one unit of time, and that they can freely allocate among themselves. The nonempty subsets of \(N\) are called \emph{coalitions}, and their set is denoted by \( \mathcal{N} \). 

\begin{definition}[Cooperative games, \textcite{von1944theory}] \leavevmode \newline
  A game is an ordered pair \((N, v)\) where 
  \begin{enumerate}
    \item \(N\) is a finite set of \emph{players}, called the \emph{grand coalition},
    \item \(v\) is a set function \( v: 2^N \to \mathbb{R} \) such that \(v(\emptyset) = 0\). 
  \end{enumerate}
\end{definition}

As previously discussed, a necessary condition for the grand coalition to form is that it must be able to distribute more utility, assuming it is transferable, among the players than any of its subcoalitions can. Denote by \(X(v)\) the affine hyperplane defined by 
\[
X(v) = \{x \in \mathbb{R}^N \mid x(N) = v(N)\}, 
\]
where \(x(S) = \sum_{i \in S} x_i\) for any coalition \(S \subseteq N\) and \(\mathbb{R}^N\) denotes the \(n\)-fold Cartesian product of \(n\) copies of \( \mathbb{R} \), one for each player. For any player \(i \in N\), the coordinate \(x_i\) represents the \emph{payment} of \(i\) by \(x\). The payment of a coalition \(S\) is simply the sum \(x(S)\) of the payment of the players. We denote by \(\mathbf{1}^i\) the \(i\)-th vector of the canonical base of \(\mathbb{R}^N\), i.e., the vector giving a payment of \(1\) to player \(i\) and \(0\) to the other players. We denote by \(\mathbf{1}^S\) the sum of vectors \(\mathbf{1}^S \coloneqq \sum_{i \in S} \mathbf{1}^i\). 

\medskip 

The set \(X(v)\) represents all the possible ways to allocate \(v(N)\) among \(n\) players. The elements of \(X(v)\) are called the \emph{preimputations} of the game \((N, v)\), and they represent the possible \emph{states} of a given game. 

\medskip 

Therefore, the question is whether there exist some preimputations such that the payment of any coalition by these preimputations is better than what \(S\) can achieve on its own. In other words, we want to know whether the following set 
\[
C(v) = \{x \in X(v) \mid x(S) \geq v(S), \forall \hspace{1pt} S \in \mathcal{N} \},
\]
called the \emph{core} of the game \((N, v)\), is nonempty. This question was solved by \textcite{bondareva1963some} and \textcite{shapley1967balanced} independently, using the concept of \emph{balanced collections}. 

\begin{definition}[Balanced collections, \textcite{bondareva1963some}, \textcite{shapley1967balanced}] \leavevmode \newline 
  A collection of coalitions \( \mathcal{B} \subseteq \mathcal{N} \) is \emph{balanced} on \(N\) if there exists a set of positive weights \( \lambda = \{\lambda_S \mid S \in \mathcal{B} \} \) such that \(\sum_{S \in \mathcal{B}} \lambda_S \mathbf{1}^S = \mathbf{1}^N\). 
  \end{definition}

A balanced collection \( \mathcal{B} \) models a social configuration of the players in \(N\). They represent how the players can be organized, respecting two natural conditions:
\begin{enumerate}
  \item The players forming a coalition spend the same time in it, 
  \item Each player is `active' for exactly one unit of time.
\end{enumerate}
The weights \( \lambda_S \) can be thought of as the fraction of the time that the players spend in the coalition \(S\). Indeed, we can reformulate the equality in the definition above by 
\[
\sum_{\substack{S \in \mathcal{B} \\ S \ni i}} \lambda_S = 1, \qquad \forall \hspace{1pt} i \in N.
\]
The total worth of the balanced collection \( \mathcal{B} \) is determined by \( \sum_{S \in \mathcal{B}} \lambda_S v(S) \). Indeed, each coalition \(S\) has a worth of \(v(S)\) and is active for \(\lambda_S\) units of time.

\begin{example}\label{ex: bal-coll}
  Consider the game defined on \( N = \{a, b, c\} \) by \( v(S) = 0 \) if \( \lvert S \rvert \in \{0, 1\} \), \(v(S) = 0.8\) if \(\lvert S \rvert = 2\) and \(v(N) = 1\). We observe that the grand coalition is the one with the largest value and that any partition of it has a worth not exceeding \(v(N)\). But the players can still be organized in a way that gives them more utility than \(v(N)\). Indeed, if each player spends half of its time with another player, and the second half of its time with the last player, during each half, the coalition secures a utility of \(0.4\), and this three times. In this case, the balanced collection is 
  \[
  \mathcal{B} = \{ \{a, b\}, \{a, c\}, \{b, c\} \} \text{ with weights being } \lambda_{\{a, b\}} = \lambda_{\{a, c\}} = \lambda_{\{b, c\}} = \frac{1}{2}
  \]
  and \(\sum_{S \in \mathcal{B}} \lambda_S v(S) = \frac{3}{2} \cdot 0.8 = 1.2 > 1\). 
\end{example}

The simple, yet rich, definition of the balanced collections made them of specific interest in other fields of mathematics, mostly in combinatorics, under the name of regular hypergraphs, or perfect fractional matching. They are also closely related to uniform hypergraphs, both the resonance and the braid hyperplane arrangements, and they form a \emph{combinatorial species of structures}~ \cite{EPTCS403.27}, combinatorial formalism created by \textcite{joyal1981theorie} allowing to use algebraic methods to solve, for example, questions about enumeration and generation of combinatorial structures. 

\medskip 

In Example~\ref{ex: bal-coll}, we cannot expect the grand coalition \(N\) to form, even if the coalition has the largest worth and no partition of it exceeds its worth. The players have a better way to be organized than forming a non-decomposable block \(N\). 

\medskip 

A game is \emph{balanced} if the balanced collection \( \{N\} \) is one with the maximal worth. In other words, a game \( (N, v) \) is balanced if, for all balanced collection \( \mathcal{B} \), we have 
\[
\sum_{S \in \mathcal{B}} \lambda_S v(S) \leq v(N). 
\]

\begin{theorem}[\textcite{bondareva1963some}, \textcite{shapley1967balanced}]\label{th: bondareva-shapley} \leavevmode \newline
The core of a game \( (N, v) \) is nonempty if and only if \( (N, v) \) is balanced. 
\end{theorem}

To put it differently, the existence of an allocation of \(v(N)\) among the players such that all coalitions receive a satisfactory payment, is equivalent to \( \{N\} \) being one of the most efficient organizations for the players in \(N\). 

\medskip 

\begin{remark}\label{remark: efficient}
The coalitions in balanced collections \( \mathcal{B} \) that satisfy \( \sum_{S \in \mathcal{B}} \lambda_S v(S) = v(N) \) are important, as they are the coalitions that cannot expect to acquire strictly more utility by cooperating. We denote their set by \( \mathcal{E}(v) \). They also play an important role in the geometry of the core, as we will see later. 
\end{remark}

Even if solutions that benefit everyone exist, they may not be accessible from the current state of the game. Assume that the players' initial endowments are the payments of an initial state represented by the preimputation \(x\). Several subcoalitions of \(N\) can be unsatisfied by this payment but may not be able to coordinate their response, because their interests are too different. Considering that the sum of all the payments must not change, an increase in the payments of some players induces a decrease in the payments of other players, leading to the inability for these coalitions to agree on a counterproposal for \(x\) and hence improve their situation simultaneously. At this stage, nothing ensures that the grand coalition \(N\) will form~\cite{laplace2024formation}. 

\medskip 

Naturally, according to the payment they receive, the coalitions may prefer a preimputation \(x\) to a preimputation \(y\) if each player in the coalition gets a strictly better payment at \(x\) than at \(y\). On the other hand, a preimputation \(x\) is affordable for a coalition \(S\) if, by working on their own, they can produce the \( \lvert S \rvert \)-dimensional vector \( x_{|S} = {(x_i)}_{i \in S} \), i.e., if \(x(S) \leq v(S)\). 

\begin{definition}[Domination, \textcite{von1944theory}]\label{def: domination} \leavevmode \newline
  Let \((N, v)\) be a game, \(x, y\) be two preimputations of \((N, v)\) and let \( S \in \mathcal{N} \) be a coalition. We say that \(x\) dominates \(y\) via \(S\) with respect to \((N, v)\), denoted \(x \hspace{4pt} \text{dom}_S^v \hspace{4pt} y\) if
  \begin{enumerate}
    \item\label{item: affordable} \(x\) is \emph{affordable} to \(S\) w.r.t. \((N, v)\), i.e, \(x(S) \leq v(S)\), 
    \item \(x\) \emph{improves} \(y\) on \(S\), i.e., \(x_i > y_i\) for all \(i \in S\). 
\end{enumerate}
\end{definition}

If the considered game \((N, v)\) is clear from the context, we simply write \(x \hspace{4pt} \text{dom}_S \hspace{4pt} y\). We say that \(x\) \emph{dominates} \(y\), denoted \(x \hspace{4pt} \text{dom} \hspace{4pt} y\), if there exists \( S \in \mathcal{N} \) such that \(x \hspace{4pt} \text{dom}_S \hspace{4pt} y\). 

\medskip

In another paper, the author~\cite{laplace2024formation} studies the steadiness of the formation of the grand coalition \( N \), by studying its internal dynamics described above. In this complementary paper, the approach is different, as we expect a third party to impose an optimal side payment to the game given its initial state. 

\medskip 

In the aforementioned paper is defined a class of polytopes, which we use in this paper, called the \emph{cooperahedra}. Intuitively, a cooperahedron is a polyhedron for which there exists a game \( (N, v) \) such that \( C(v) \) is its topological closure. More formally, we have the following definition. Let \( \langle \cdot, \cdot \rangle \) denote the usual scalar product of \( \mathbb{R}^N \). 

\begin{definition}[Cooperahedra, \cite{laplace2024formation}] \leavevmode \newline 
A \emph{cooperahedron} is a polyhedron \( P \subseteq \mathbb{R}^N \) of the form
\[
P = \{ x \in \mathbb{R}^N \mid x(N) = b^0 \text{ and } \langle z^i, x \rangle \diamond_i b^i, \forall i \in I_P \},
\]
where \( I_P \) is finite, and for all \( i \in I_P \), we have \( \diamond_i \in \{\geq, >, <, \leq \} \) and \( z^i \in \{0, 1\}^N \). 
\end{definition}

For convenience, we denote by \( I^\prec_P \) and \( I^\succ_P \) the sets of indices defined by 
\[
I^\prec_P = \{ i \in I_P \mid \diamond_i \in \{<, \leq \} \} \quad \text{and} \quad I^\succ_P = \{ i \in I_P \mid \diamond_i \in \{\geq, >\} \}. 
\]
For each \( i \in I_P \), we can associate a coalition \( S_i \) defined by 
\[
S_i = \begin{cases} 
\{j \in N \mid z^i_j = 1 \}, & \quad \text{if } i \in I^\succ_P, \\
\{j \in N \mid z^i_j = 0 \}, & \quad \text{if } i \in I^\prec_P.
\end{cases}
\]
Hence, with each cooperahedron, we can define a game \( (\mathcal{F}_P, v_P) \) on a set system \( \mathcal{F}_P \), which can be a proper subcollection of \( \mathcal{N} \), by 
\[
\mathcal{F}_P = \{S_i \mid i \in I\} \qquad \text{and} \qquad v_P(S_i) = \begin{cases}
b_i, & \quad \text{if } i \in I^\succ_P, \\
b^0 - b^i, & \quad \text{if } i \in I^\prec_P. 
\end{cases}
\]
The next result is the reason why we use the cooperahedron in this paper. Let \( I^*_P \) be the subset of \( I_P \) defined by \( I^*_P \coloneqq \{ i \in I_P \mid \diamond_i \in \{>, <\} \} \). 

\begin{theorem}[\cite{laplace2024formation}]\label{th: main} \leavevmode \newline
A cooperahedron \( P \) is nonempty if and only if \( (\mathcal{F}_P, v_P) \) is balanced and 
\[
\mathcal{E} \left( \mathcal{F}_P, v_P \right) \cap \{ S_i \mid i \in I^*_P \} = \emptyset. 
\]
\end{theorem}

The class of cooperahedra includes several well-known classes of polytopes and polyhedra, the most famous of them being the \emph{deformed permutohedra} \cite{postnikov2009permutohedra}, related to convex games. For more details about cooperahedra, see \cite{laplace2024formation}. 
 
\section{The geometry of a cooperative game}\label{section: geometry}

This section focuses on the affine subspace \( X(v) \subseteq \mathbb{R}^N \). The linear subspace of \( \mathbb{R}^N \) parallel to \( X(v) \) is denoted by \( \Sigma \), and its elements are called the \emph{side payments}. A side payment represents a translation in the space of preimputations, which is a redistribution of utility between players, without changing the total amount acquired by the players. Indeed, if \( X(v) \) is defined by \( X(v) = \{ x \in \mathbb{R}^N \mid x(N) = v(N) \} \), then \( \Sigma = \{ x \in \mathbb{R}^N \mid x(N) = 0 \} \) is a vector space. Also, notice that \( \Sigma \) does not depend on a given game, but only on the number of players, which defines its dimension. 

\medskip 

Because \( N \) is just an arbitrary coalition, we want to have the same subspaces in hands for smaller coalitions. Let \( S \in \mathcal{N} \) be a coalition. We denote by \( A_S(v) \) the set of preimputations for which \( S\) exactly gets its value, i.e., 
\[
A_S(v) \coloneqq \{ x \in X(v) \mid x(S) = v(S) \}. 
\]
The linear subspace of \( \mathbb{R}^N \) parallel to \( A_S(v) \) is denoted by \( H_S \), is a linear subspace of \( \Sigma \) as well, and is defined by \( H_S = \{\sigma \in \Sigma \mid \sigma(S) = 0 \} \). These two sets of affine subspaces define \emph{hyperplane arrangements}, respectively in \( X(v) \) and \( \Sigma \). The hyperplane arrangement formed by \( \{H_S \mid S \in \mathcal{N} \} \) is known as the \emph{resonance arrangement}, and does not depend on the game under consideration. However, the hyperplane arrangement formed by \( \{ A_S(v) \mid S \in \mathcal{N} \} \) is a distortion of the resonance arrangement that depends on the game. Hence, its \emph{chambers}, i.e., the connected components of the complement of \( \{ A_S(v) \mid S \in \mathcal{N} \} \), summarize the game by containing the right amount of information about it, especially about the core, and domination between preimputations. 

\begin{figure}[ht]
\begin{center}
\begin{subfigure}{0.49\textwidth}
\begin{center}
\begin{tikzpicture}[scale=0.25]
\draw[cyan] (-10, 0) node[left] {\footnotesize $A_{\{a, b\}}$} -- (10, 0) node[right] {\footnotesize $A_{\{c\}}$};
\draw[purple] (-5.774, -10) node[left] {\footnotesize $A_{\{a\}}$} -- (5.774, 10) node[right] {\footnotesize $A_{\{b,  c\}}$};
\draw[orange] (5.774, -10) node[right] {\footnotesize $A_{\{b\}}$} -- (-5.774, 10) node[left] {\footnotesize $A_{\{a, c\}}$};


\end{tikzpicture}
\caption{The resonance arrangement.}
\end{center}
\end{subfigure}
\begin{subfigure}{0.49\textwidth}
\begin{center}
\begin{tikzpicture}[scale=0.25]
\fill[cyan!10] (-10, 0) -- (10, 0) -- (10, -2.578) -- (-10, -2.578) -- cycle;
\fill[purple!10] (-7.706, -10) -- (3.841, 10) -- (7.788, 10) -- (-3.760, -10) -- cycle;
\fill[orange!10] (5.774, -10) -- (9.233, -10) -- (-2.314, 10) -- (-5.774, 10) -- cycle;

\fill[gray!10] (0.764, 4.669 ) -- (2.737, 1.251) -- (2.014, 0) -- (0, 0) -- (-0.966, 1.673) -- cycle;
\fill[gray!10] (4.948, -2.578) -- (1.488, -2.578) -- (1.007, 1.744) -- (2.014, 0) -- (3.459, 0) -- cycle;
\fill[gray!10] (-3.420, -2.578) -- (-1.932, 0) -- (0, 0) -- (1.007, 1.744) -- (0.526, -2.578) -- cycle;

\fill[gray!50] (0, 0) -- (2.014, 0) -- (1.007, -1.744) -- cycle;

\draw[cyan] (-10, 0) node[left] {\footnotesize $A_{\{a, b\}}$} -- (10, 0);
\draw[cyan] (-10, -2.578) -- (10, -2.578) node[right] {\footnotesize $A_{\{c\}}$};
\draw[purple] (-7.706, -10) node[left] {\footnotesize $A_{\{a\}}$} -- (3.841, 10); 
\draw[purple] (-3.760, -10) -- (7.788, 10) node[right] {\footnotesize $A_{\{b, c\}}$};
\draw[orange] (5.774, -10) -- (-5.774, 10) node[left] {\footnotesize $A_{\{a, c\}}$};
\draw[orange] (-2.314, 10) -- (9.233, -10) node[right] {\footnotesize $A_{\{b\}}$};


\end{tikzpicture}
\caption{Arrangement \( \{ A_S(v) \mid S \in \mathcal{N} \} \).}
\end{center}
\end{subfigure}
\caption{Distortion of the resonance arrangement by a game $(N, v)$. The core $C(v)$ is the darkest area.}
\label{fig: slab}
\end{center}
\end{figure}
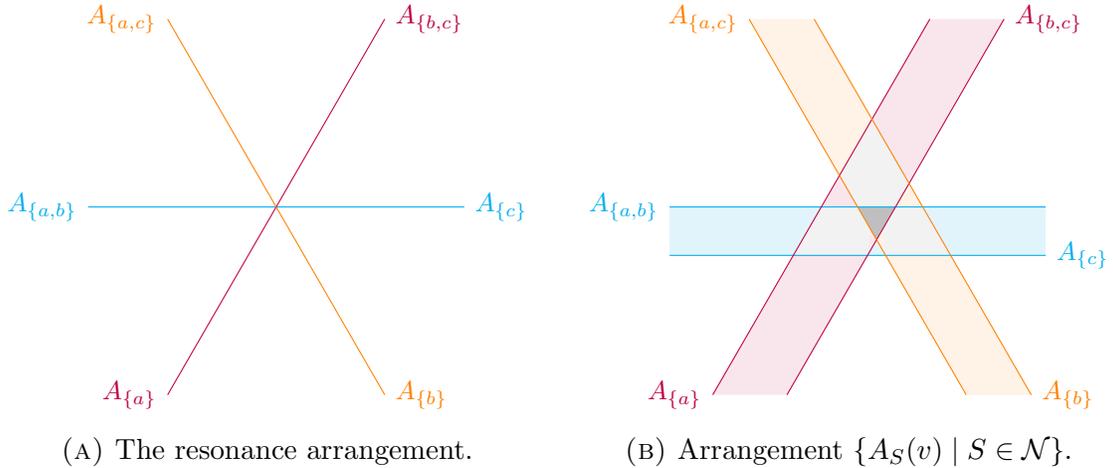

The coalitions in \( \mathcal{E}(v) \) (see Remark~\ref{remark: efficient}) define important affine subspaces. 

\begin{proposition}[\cite{laplace2024formation}]\label{prop: efficient_hyperplane} \leavevmode \newline
Let \( (N, v) \) be balanced. A coalition \( S \) belongs to \( \mathcal{E}(v) \) if and only if \( C(v) \subseteq A_S(v) \). 
\end{proposition}

By definition, a preimputation \( x \) can be dominated via a coalition \( S \) only if the coalition can \emph{improve upon} \( x \) by leaving the grand coalition, i.e., \( x(S) < v(S) \). In the sequel, we denote by \( e_S(x) \) the additional quantity of money that coalition \( S \) can acquire by itself, i.e., \( e_S(x) \coloneqq v(S) - x(S) \), that we call the \emph{excess of \( S \) at \( x \)}. 

\begin{definition}[Feasible collections, \textcite{grabisch2021characterization}]\label{def: feasible-collections} \leavevmode \newline
Let \( \mathcal{Q} \subseteq \mathcal{N} \) be a collection of coalitions. We denote by \( X_\mathcal{Q}(v) \) the set of preimputations upon which a coalition can improve if and only if it belongs to \( \mathcal{Q} \), i.e., 
\[
X_\mathcal{Q}(v) \coloneqq \{x \in X(v) \mid x(S) < v(S) \text{ if and only if } S \in \mathcal{Q} \}. 
\]
We say that the collection \( \mathcal{Q} \) is \emph{feasible} if its associated \emph{region} \( X_\mathcal{Q}(v) \) is nonempty.
\end{definition} 

\begin{remark}
The region \( X_\mathcal{Q}(v) \) is a cooperahedron. 
\end{remark}

Note that the core is the region associated with the empty collection of coalitions \( \mathcal{Q} = \emptyset \), hence the nonempty regions form a partition of \( X(v) \). Also, note that the notion of chambers and regions are quite similar. However, not every region is a chamber, as they can have a lower dimension than the ones of the chambers. 

\medskip 

Throughout this paper, we equip the Euclidean vector space \( \mathbb{R}^N \), and therefore \( X(v) \), with the usual scalar product, denoted by \( \langle \cdot, \cdot \rangle \) and the associated norm, denoted by \( \lVert \cdot \rVert \) defined, for all \( x, y \in \mathbb{R}^N \), by 
\[
\langle x, y \rangle = \sum_{i \in N} x_i y_i, \qquad \text{and} \qquad \lVert x \rVert = \sqrt{ \langle x, y \rangle }. 
\]
A vector \( x \in \mathbb{R}^N \) is said to be a \emph{normal} of a subset \( Y \subseteq \mathbb{R}^N \) if, for every element \( y \in Y' \) with \( Y' \) being the linear subspace of \( \mathbb{R}^N \) parallel to \( Y \), we have \( \langle x, y \rangle = 0 \). The normal vectors are very useful to handle hyperplanes, as they greatly simplify formulae. 

\begin{proposition}
Let \( S \) be a coalition. The vector \( \eta^S \), defined by 
\[
\eta^S = \mathbf{1}^S - \frac{\lvert S \rvert}{\lvert N \rvert} \mathbf{1}^N, 
\]
is a side payment and a normal of \( A_S(v) \). 
\end{proposition}

\begin{proof}
For any coalition \( S \), we have 
\[
\eta^S(N) = \mathbf{1}^S (N) - \frac{\lvert S \rvert}{\lvert N \rvert} \mathbf{1}^N = \lvert S \rvert - \lvert S \rvert = 0, 
\]
hence \( \eta^S \) is a side payment. Let \( x \) and \( y \) be two elements of \( A_S(v) \). We have 
\[
\langle \eta^S, x - y \rangle = \langle \eta^S, x \rangle - \langle \eta^S, y \rangle = x(S) - \frac{\lvert S \rvert}{\lvert N \rvert} x(N) - y(S) + \frac{\lvert S \rvert}{\lvert N \rvert} y(N). 
\]
As \( x(S) = y(S) = v(S) \) and \( x(N) = y(N) = v(N) \), we have \( \langle \eta^S, x - y \rangle = 0 \). 
\end{proof}

We extend the definition of these vectors by setting \( \eta^\emptyset = 0_{\mathbb{R}^N} \), which coincides with \( \eta^N \). The set \( \{ \eta^S \mid S \subseteq N \} \) does not depend on the game we study, and the coalition function of the game only defines the position of \( A_S(v) \) along the line generated by \( \eta^S \). 

\medskip 

Now that we have a precise and handy characterization of the affine subspaces \( A_S(v) \), we want to define projectors onto them. 

\begin{theorem}[Hilbert Projection Theorem]\label{th: hilbert-proj} \leavevmode \newline 
For every \( x \in \mathbb{R}^N \) and every nonempty closed convex \( K \subseteq \mathbb{R}^N \), there exists a unique element \( z \in K \) for which \( \lVert x - z \rVert \) is equal to \( \inf_{y \in K} \lVert x - y \rVert \). 
\end{theorem}

The element \( z \in K \) is called the \emph{projection} of \( x \) onto \( K \). We denote by \( \pi_K: \mathbb{R}^N \to K \) the map that assigns to \( x \in \mathbb{R}^N \) its projection \( \pi_k(x) \) onto \( K \). The map \( \pi_K \) is called the \emph{projector} onto \( K \). In order to compute projection, we use the following characterization in the formulation of \textcite{bauschke2011convex}. 

\begin{theorem}[Projection Theorem] \leavevmode \newline 
Let \( K \) be a nonempty closed convex subset of \( \mathbb{R}^N \). Then, for all \( x \) and \( y \) in \( \mathbb{R}^N \), 
\[
y = \pi_K(x) \quad \text{if and only if} \quad \Big[ y \in K \quad \text{and for all } z \in K, \langle z - y, x - y \rangle \leq 0 \Big].
\]
\end{theorem}

We now build our projectors onto each affine subspace \( A_S(v) \). 

\begin{proposition} \label{prop: first-proj}
Let \( S \in \mathcal{N} \setminus \{ N \} \) be a coalition and \( x \) be a preimputation. Then, 
\[
\pi_{A_S(v)} = x + \gamma_S(x) \eta^S, 
\]
where \( \gamma_S(x) \coloneqq e_S(x) / \lVert \eta^S \rVert^2 \). Moreover, if \( x(S) < v(S) \), we have \( \pi_{A_S(v)} \hspace{4pt} \text{dom}_S \hspace{4pt} x \). 
\end{proposition}

\begin{proof}
First, notice that \( \lVert \eta^S \rVert^2 = \langle \eta^S, \eta^S \rangle = \eta^S(S) \). Following the Projection Theorem, we first prove that \( \pi_{A_S(v)} \) belongs to \( A_S(v) \). 
\[
\left( x + \gamma_S(x) \eta^S \right)(S) = x(S) + \gamma_S(x) \eta^S(S) = x(S) + \frac{e_S(x)}{\lVert \eta^S \rVert^2} \lVert \eta^S \rVert^2 = x(S) + e_S(x) = v(S). 
\]
Let \( z \) be an element of \( A_S(v) \). We have 
\[ \begin{aligned} 
\big\langle z - (x + \gamma_S(x) \eta^S, \hspace{2pt} x - (x + \gamma_S(x) \eta^S) \big\rangle & = \langle z - x - \gamma_S(x) \eta^S, \hspace{2pt} - \gamma_S(x) \eta^S \rangle \\
& = - \gamma_S(x) \left( \langle z, \eta^S \rangle - \langle x, \eta^S \rangle - \gamma_S(x) \langle \eta^S, \eta^S \rangle \right) \\
& = - \gamma_S(x) \left( z(S) - x(S) - \gamma_S(x) \lVert \eta^S \rVert^2 \right) \\
& = - \gamma_S(x) \left( e_S(x) - e_S(x) \right) \\
& = 0.
\end{aligned} \]
We use the Projection Theorem to conclude the construction of the projector. To end this proof, notice that, if \( e_S(x) > 0 \), for all \( i \in S \), we have \( \eta^S_i = 1 - \lvert S \rvert / n > 0 \). 
\end{proof}

\begin{figure}[ht]
\begin{center}
\begin{tikzpicture}[scale=0.2]
\draw (1.774, 10) node[above left]{$A_S(v)$} -- (11.774, -10);

\filldraw[black] (-16,-6) circle (5pt) node[above left]{$x$};
\filldraw[black] (4.6192, 4.3096) circle (5pt) node[right]{$\pi_{A_S(v)}(x)$};

\draw[dashed] (-16, -6) -- (4.6192, 4.3096);
\draw (2.7178, 4.9434) -- (3.9854, 5.5772);
\draw (2.7178, 4.9434) -- (3.3516, 3.6758);

\draw[{Stealth[scale=1.5]}-{Stealth[scale=1.5]}] (-15.3763, -7.2474) -- (5.2429, 3.0622); 

\draw[-{Stealth[scale=1.2]}] (9.8192, -6.0904) -- (14, -4) node[below]{\small{$\eta^S$}};

\path (-3.3, -2.7) node[below]{\small{$\gamma_S(x) = \dfrac{e_S(x)}{\lVert \eta^S \rVert}$}};
\end{tikzpicture}
\caption{Projection of $x$ onto $A_S(v)$.}
\label{fig: proj-1}
\end{center}
\end{figure}
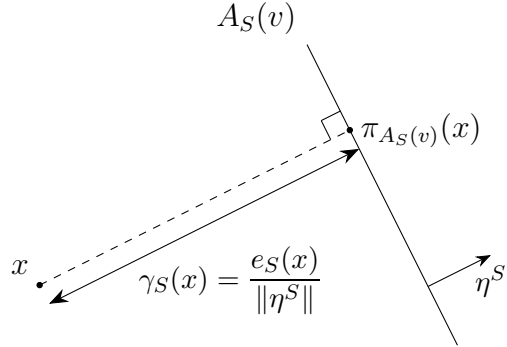

If a coalition \( S \) can improve upon a preimputation \( x \in X(v) \), the projection \( \pi_{A_S(v)}(x) \) shares the excess of \( S \) at \( x \) equally among the players in \( S \), and therefore the projection dominates \( x \) via \( S \). This fact motivates the use of these projectors in this work to study the domination relations between payment vectors and, more specifically, between preimputations and core elements. Moreover, \( \gamma_S(x) \eta^S \) is the shortest side payment able to map the preimputation \( x \) onto a preimputation that benefits the coalition \( S \) (see Figure \ref{fig: proj-1}). It is the solution to the problem consisting of satisfying coalition \( S \) with the shortest (with respect to the Euclidean distance) side payment. 

\medskip

The idea of projections between preimputations was already existing in the transfer scheme defined by \textcite{cesco1998convergent}. A transfer scheme, first defined and used by \textcite{stearns1968convergent}, is a sequence of preimputations defined by a sequence of transfers, i.e., translation by a side payment. The sequence of preimputations defined by Cesco is defined recursively by projecting onto \( A_S(v) \) the last preimputation produced. His projector is defined, for every preimputation \( x \), by 
\[
\pi_{A_S(v)}(x) = x + e_S(x) \left( \frac{\mathbf{1}^S}{\lvert S \rvert} - \frac{\mathbf{1}^{N \setminus S}}{\lvert N \setminus S \rvert} \right). 
\]
Cesco's projectors are identical to ours, but the choice of normals is different. He chose to have a normal of a different norm to have a formulation depending explicitly on the excess \( e_S(x) \). Our choice of normals is motivated by the following theorem. 

\begin{theorem} \label{th: modular_normal}
Let \( S \) and \( T \) be coalitions. We have 
\[
\eta^S + \eta^T = \eta^{S \cup T} + \eta^{S \cap T} \qquad \text{and} \qquad \langle \eta^S, \eta^T \rangle = \eta^S(T). 
\]
\end{theorem} 

\proof Notice that if \( S \) and \( T \) are disjoint, we have 
\[
\eta^S + \eta^T = \mathbf{1}^S - \frac{\lvert S \rvert}{\lvert N \rvert} \mathbf{1}^N + \mathbf{1}^T - \frac{\lvert T \rvert}{\lvert N \rvert} \mathbf{1}^N = \mathbf{1}^{S \cup T} - \frac{\lvert S \rvert + \lvert T \rvert}{\lvert N \rvert} \mathbf{1}^N = \eta^{S \cup T}. 
\]
Otherwise we decompose \( \eta^S = \eta^{S \setminus T} + \eta^{S \cap T} \), similarly for \( \eta^T \), and we have
\[
\eta^S + \eta^T = \eta^{S \setminus T} + \eta^{S \cap T} + \eta^{T \setminus S} + \eta^{S \cap T} = \eta^{S \cup T} + \eta^{S \cap T}. 
\]
For the second property, the fact that \( \eta^S \) is a side payment leads to
\[
\langle \eta^S, \eta^T \rangle = \langle \eta^S, \mathbf{1}^T \rangle - \frac{\lvert T \rvert}{\lvert N \rvert} \langle \eta^S, \mathbf{1}^N \rangle = \eta^S(T) - \frac{\lvert T \rvert}{\lvert N \rvert} \eta^S(N) = \eta^S(T).
\] 
\endproof

From this theorem, we can derive several properties on the set \( \{ \eta^S \mid S \in \mathcal{N} \} \), in particular, the repartition of the \( \eta^S \) is balanced around the origin of \( \Sigma \).  

\begin{theorem} \label{th: balancedness_zero}
A collection of coalitions \( \mathcal{C} \subseteq \mathcal{N} \) is balanced if and only if there exists a set of positive numbers \( \{\theta_S \mid S \in \mathcal{C}\} \) such that 
\[
\sum_{S \in \mathcal{C}} \theta_S \eta^S = 0_{\mathbb{R}^N}.
\]
\end{theorem}

\proof
Let \( \mathcal{C} \) be a balanced collection. Then, there exists a set of positive numbers \( \{ \lambda_S \mid S \in \mathcal{C} \} \) such that 
\[
\sum_{S \in \mathcal{C}} \lambda_S \mathbf{1}^S = \mathbf{1}^N. 
\]
Taking the sum of the coordinates of both sides gives that following identity:
\[
\sum_{S \in \mathcal{C}} \lambda_S \lvert S \rvert = \lvert N \rvert. 
\]
Then we have 
\[
\sum_{S \in \mathcal{C}} \lambda_S \eta^S = \sum_{S \in \mathcal{C}} \lambda_S \mathbf{1}^S - \frac{1}{\lvert N \rvert} \sum_{S \in \mathcal{C}} \lambda_S \lvert S \rvert \mathbf{1}^N = \mathbf{1}^N - \frac{1}{\lvert N \rvert} \lvert N \rvert \mathbf{1}^N = 0_{\mathbb{R}^N}. 
\]
Let assume now assume that there exists positive numbers \( \{ \theta_S \mid S \in \mathcal{C} \} \) such that 
\[
\sum_{S \in \mathcal{C}} \theta_S \eta^S = 0_{\mathbb{R}^N},
\]
and set \( \alpha = \frac{1}{\lvert N \rvert} \sum_{S \in \mathcal{C}} \theta_S \lvert S \rvert > 0 \). By developing the vectors \( \eta^S \), we get 
\[
\sum_{S \in \mathcal{C}} \theta_S \eta^S = \sum_{S \in \mathcal{C}} \theta_S \left( \mathbf{1}^S - \frac{\lvert S \rvert}{\lvert N \rvert} \mathbf{1}^S \right) = \sum_{S \in \mathcal{C}} \theta_S \mathbf{1}^S - \mathbf{1}^N \frac{1}{\lvert N \rvert} \sum_{S \in \mathcal{C}} \theta_S \lvert S \rvert = \sum_{S \in \mathcal{C}} \theta_S \mathbf{1}^S - \alpha \mathbf{1}^N. 
\]
By hypothesis, we then have that \( \sum_{S \in \mathcal{C}} \theta_S \mathbf{1}^S - \alpha \mathbf{1}^N = 0 \), and by dividing both sides by \( \alpha \) and subtracting by \( \mathbf{1}^N \), we get 
\[
\sum_{S \in \mathcal{C}} \alpha^{-1} \theta_S \mathbf{1}^S = \mathbf{1}^N, 
\]
so \( \mathcal{C} \) is a balanced collection with weights \( \{\alpha^{-1} \theta_S \mid S \in \mathcal{C}\} \). 
\endproof

We can interpret the vectors \( \eta^S \) as the direction along we can translate a preimputation such that the increase of the payment of \( S \) is maximized. In other words, 
\[
\frac{1}{\lVert \eta^S \rVert} \eta^S = \argmax_{x \in B_\Sigma (0, 1)} x(S),
\]
with \( B_\Sigma (0, 1) \coloneqq \{x \in \Sigma \mid \lVert x \rVert = 1\} \) is the unit ball of the space of side payments. Hence, each normal vector represents the interest of its corresponding coalition, as it is the gradient of the linear form of its payment function. A balanced collection is therefore characterized by the fact that the interests of the coalitions included in it are perfectly incompatible. Because a convex combination of these vectors is zero, any increase in the direction profiting one coalition necessarily decreases the payment of another coalition in this coalition. This interpretation of the balanced collections resembles a Pareto optimality of the social structure formed by the players. 

\medskip 

The opposite property characterizes the \emph{unbalanced collections}. A collection \( \mathcal{C} \subseteq \mathcal{N} \) is called \emph{unbalanced} if it does not contain a balanced collection. In particular, it is not itself a balanced collection. It is important to notice that unbalancedness is not equivalent to the negation of balancedness. The collection \( \{ \{a, b\}, \{a, c\}, \{b, c\}, \{a\} \} \) is not balanced on \( \{a, b, c\} \), but contains \( \{ \{a, b\}, \{a, c\}, \{b, c\} \} \), which is balanced. 

\medskip 

These collections that are not balanced but contain a balanced collection are called \emph{weakly balanced} by \textcite{maschler1971kernel}. This is motivated by two aspects: first, by removing some coalitions, a weakly balanced collection can become balanced. Secondly, relaxing the positivity condition for a non-negativity condition on the weights in the definition of balancedness leads to weak balancedness.  

\medskip 

A unbalanced collection which is maximal with respect to inclusion is called a \emph{maximal unbalanced collection}. Contrarily to the balanced collections, an unbalanced collection represents a set of coalitions with common interests. 

\begin{proposition}[Billera, Moore, Moraites, Wang and Williams \cite{billera2012maximal}] \label{prop: unbalanced} \leavevmode \newline
A collection \( \mathcal{S} \subseteq \mathcal{N} \) is unbalanced if and only if there exists a side payment \( \sigma \in \Sigma \) such that, for all \( S \in \mathcal{S} \), we have \( \sigma(S) > 0 \). 
\end{proposition}

\begin{example}
These two collections are maximal unbalanced collections with \( \sigma \) a side payment illustrating the previous result: 
\begin{itemize}
\item for \( n = 3 \): \( \{ \{a, b\}, \{a, c\}, \{a\} \} \) and, as a side payment, \( \sigma = (2, -1, -1) \);
\item for \( n = 4 \): \( \{ \{a\}, \{a, b\}, \{a, c\}, \{a, d\}, \{a, b, c\}, \{a, b, d\}, \{a, c, d\} \} \) and, as a side payment, \( \sigma = (3, -1, -1, -1) \). 
\end{itemize}
The side payment \( \sigma \) represents a direction which improves the payment of all coalitions included in the unbalanced collection. 
\end{example}

To conclude this section, let us have a look at an interesting corollary of Theorem~\ref{th: modular_normal} about the modularity of the normal vectors. For any coalition \( S \) and \( T \) both distinct from \( N \), if \( \lvert N \rvert \) is a prime number, \( \eta^S \) and \( \eta^T \) cannot be orthogonal, because
\[
\langle \eta^S, \eta^T \rangle = \eta^S(T) = \mathbf{1}^S(T) - \frac{\lvert S \rvert}{\lvert N \rvert} \mathbf{1}^N(T) = \lvert S \cap T \rvert - \frac{\lvert S \rvert \cdot \lvert T \rvert}{\lvert N \rvert},
\] 
and $\lvert S \rvert \cdot \lvert T \rvert$ cannot be a multiple of $\lvert N \rvert$. Then, no hyperplane in the resonance arrangement can be orthogonal. From a game-theoretic point of view, it means that when the number of players is a prime number, any optimal increase in the payment of a coalition, i.e., the current allocation is shifted in the direction of \( \eta^S \) for a given \( S \in \mathcal{N} \), induces a nonzero evolution of the payment of every coalition. 

\section{Projection onto the intersection of \( A_S(v) \) and \( A_T(v) \)}


Thanks to Proposition~\ref{prop: first-proj}, we know that for a preimputation \( x \in X(v) \), if \( x(S) < v(S) \), we have \( \pi_{A_S(v)} (x) \hspace{4pt} \text{dom}_S \hspace{4pt} x \). But we have no control on where the projection lies. For any preimputation \( x \), the relative position of \( \pi_{A_S(v)}(x) \) with respect to \( A_T(v) \) depends on the excesses \( e_S(x) \) and \( e_T(x) \), but also on the scalar product \( \eta^S(T) = \langle \eta^S, \eta^T \rangle \).

\medskip 


The value and the sign of \( \langle \eta^S, \eta^T \rangle \) indicate how correlated the evolution of the payments of coalitions \( S \) and \( T \) are, i.e., how much their interests overlap. Adding the side payment \( \eta^S \) to a preimputation \( x \) necessarily benefits coalition \( S \), and the value of \( \langle \eta^S, \eta^T \rangle \) indicates how much the translation by \( \eta^S \) benefits coalition \( T \). Indeed, 
\[
\left( x + \eta^S \right)(T) = x(T) + \eta^S(T) = x(T) + \langle \eta^S, \eta^T \rangle. 
\]
To summarize all this information, we define 
\[
\chi_S(T, x) \coloneqq e_S(x) \langle \eta^S, \eta^T \rangle - e_T(x) \lVert \eta^S \rVert^2.
\]

\begin{proposition} \label{prop: chi}
Let \( S \) and \( T \) be two coalitions, and let \( x \in X(v) \). Then
\[
\pi_{A_S(v)}(x) \in A^\geq_T \quad \text{if and only if} \quad \chi_S(T, x) \geq 0, 
\]
with \( A^\geq_T = \{ x \in X(v) \mid x(T) \geq v(T) \} \). 
\end{proposition}

\proof
First, we study the excess of \( T \) at the projection onto \( A_S(v) \):
\[ \begin{aligned}
e_T(\pi_{A_S(v)}(x)) & = v(T) - x(T) - \gamma_S(x) \eta^S(T) \\
& = e_T(x) - \frac{e_S(x)}{\lVert \eta^S \rVert^2} \langle \eta^S, \eta^T \rangle \\
& = \frac{1}{\lVert \eta^S \rVert^2} \left( e_T(x) \lVert \eta^S \rVert^2 - e_S(x) \langle \eta^S, \eta^T \rangle \right) = \frac{-1}{\lVert \eta^S \rVert^2} \chi_S(T, x). 
\end{aligned} \]
The projection lies into \( A^\geq_T \) if and only if \( e_T(\pi_{A_S(v)}(x)) \) is non-positive and therefore if and only if \(\chi_S(T, x)\) is non-negative. 
\endproof

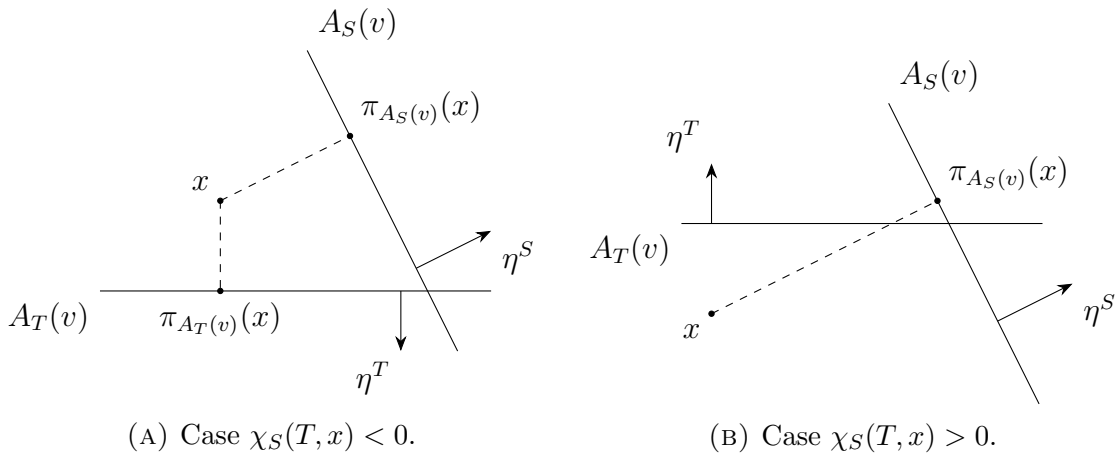
\begin{figure}[ht]
\begin{center}
\begin{subfigure}{0.49\textwidth}
\begin{center}
\begin{tikzpicture}[scale=0.2]
\draw (1.774, 10) node[above right]{\( A_S(v) \)} -- (11.774, -10);
\draw (-12, -6) node[below left]{\( A_T(v) \)} -- (14, -6);

\filldraw[black] (-4,0) circle (5pt) node[above left]{\( x \)};
\filldraw[black] (4.6192, 4.3096) circle (5pt) node[above right]{\( \pi_{A_S(v)}(x) \)};
\filldraw[black] (-4, -6) circle (5pt) node[below]{\( \pi_{A_T(v)}(x) \)};

\draw[dashed] (-4, -0) -- (4.6192, 4.3096);
\draw[dashed] (-4, 0) -- (-4, -6);

\draw[-{Stealth[scale=1.2]}] (9.0192, -4.4904) -- (14, -2) node[below right]{\( \eta^S \)};
\draw[-{Stealth[scale=1.2]}] (8, -6) -- (8, -10) node[below left]{\( \eta^T \)};
\end{tikzpicture}
\caption{Case \( \chi_S(T, x) < 0 \).}
\end{center}
\end{subfigure}
\begin{subfigure}{0.49\textwidth}
\begin{center}
\begin{tikzpicture}[scale=0.2]
\draw (1.774, 10) node[above right]{\( A_S(v) \)} -- (11.774, -10);
\draw (-12, 2) node[below left]{\( A_T(v) \)} -- (12, 2);

\filldraw[black] (-10, -4) circle (5pt) node[below left]{\( x \)};
\filldraw[black] (5.0192, 3.5096) circle (5pt) node[above right]{\( \pi_{A_S(v)}(x) \)};

\draw[dashed] (-10, -4) -- (5.0192, 3.5096);

\draw[-{Stealth[scale=1.2]}] (9.0192, -4.4904) -- (14, -2) node[below right]{\( \eta^S \)};
\draw[-{Stealth[scale=1.2]}] (-10, 2) -- (-10, 6) node[above left]{\( \eta^T \)};
\end{tikzpicture}
\caption{Case \( \chi_S(T, x) > 0 \).}
\label{subfig: chi_positive}
\end{center}
\end{subfigure}
\caption{Illustrative example of Proposition \ref{prop: chi}.}
\label{fig: chi}
\end{center}
\end{figure}

There are intuitive implications for this result. If \( \langle \eta^S, \eta^T \rangle \geq 0 \), increasing the payment of coalition \( S \) by translating a given preimputation along \( \eta^S \) also increases the payment of coalition \( T \). Therefore, assuming that \( e_S(x) \) is sufficiently large, the side payment \( \sigma \) between \( x \) and \( \pi_{A_S(v)}(x) \) increases the payment of \( S \) as well as the payment of \( T \), and we have \( \pi_{A_S(v)}(x)(T) = x(T) + \sigma(T) \geq v(T) \), as we can see on Figure~\eqref{subfig: chi_positive}. 

\medskip

For most of the preimputations, projecting onto an affine subspace \( A_S(v) \) is not sufficient to project onto the core. To do so, we need to be able to define projectors on any intersection of affine subspaces, to simultaneously improve the payment of all the coalitions that can improve upon the considered preimputation. 

\section{Projections onto arbitrary intersections of affine subspaces \(A_S(v)\)}\label{sec: projection-arbitrary-intersection}

In the previous section, we described a characterization which allows us to know whether there exists a side payment simultaneously improving the payment of a collection of coalitions, and if this side payment could be supported by a domination process involving one of the coalitions in the collection. In this section, we are looking for the smallest (according to the Euclidean norm) side payment between an affine subspace of \( X(v) \) and a given preimputation. 

\medskip

Let us consider a preimputation \( x \in X(v) \), an element \( y \) of a given affine subspace \( V \subseteq X(v) \), and the side payment \( \sigma \) such that \( x + \sigma = y \). Finding the side payment with the minimal norm is equivalent to, for the given preimputation \( x \), solving this minimization problem \( \min_{y \in V} \lVert y - x \rVert \). By definition of the Euclidean norm, we have 
\[
\argmin_{y \in V} \; \lVert y - x \rVert = \argmin_{y \in V} \sqrt{ \sum_{i \in N} \left( y_i - x_i \right)^2} = \argmin_{y \in V} \sum_{i \in N} \left( y_i - x_i \right)^2. 
\]
By the Hilbert Projection Theorem (Theorem \ref{th: hilbert-proj}), the element \( y \) is the (orthogonal) projection of \( x \) onto \( V \). The terms \( (y_i - x_i) \) in the formula above are the adjustments of the individual players' payments induced by the side payment \( \sigma = y - x \). Then, \( \sigma \) is the side payment reorganizing the payments of individual players which minimize the sum of the square of the adjustments, while preserving the total sum of payments. Moreover, by the Hilbert Projection Theorem, it is the only side payment with these properties. In order to identify these side payments, we study projectors on \( X(v) \). 

\begin{definition}
A \emph{projector} \( P \) is an idempotent linear map, i.e., \( P \circ P = P \).
\end{definition}

The image of \( P: \mathbb{R}^N \to \mathbb{R}^N \), denoted by \( \text{im}(P) \) and called the \emph{column space} in the context of linear algebra, is the subspace of \( \mathbb{R}^N \) spanned by the column of \( P \). Then, \( P \) sends any element \( x \) of the space \( \mathbb{R}^N \) onto an element \( z \in \text{im}(P) \), which is the closest element of \( \text{im}(P) \) from \( x \) according to the Hilbert Projection Theorem.

\begin{remark}[Four fundamental subspaces] \label{remark: fundamental-subspaces}
Each matrix \( A \) a size \( (n \times k) \) defines four subspaces, called the \emph{four fundamental subspaces} of \( A \), which are 
\begin{itemize}[itemsep = 5pt]
\item the \emph{column space} of \( A \): \( {\rm im}(A) \coloneqq \{Ay \mid y \in \mathbb{R}^k\} \subseteq \mathbb{R}^n \), 
\item the \emph{row space} of \( A \): \( {\rm im}(A^\top) \coloneqq \{A^\top x \mid x \in \mathbb{R}^n\} \subseteq \mathbb{R}^k \), 
\item the \emph{kernel} of \( A \): \( {\rm ker}(A) \coloneqq \{y \in \mathbb{R}^k \mid Ay = 0\} \), 
\item the \emph{left kernel} of \( A \): \( {\rm ker}(A^\top) \coloneqq \{ x \in \mathbb{R}^n \mid A^\top x = 0\} \). 
\end{itemize}
These four subspaces are related in the following way (see \parencite[][Fact 2.25]{bauschke2011convex}):
\[
\mathbb{R}^n = {\rm im}(A) \oplus \ker \left( A^\top \right) \qquad \text{and} \qquad \mathbb{R}^k = {\rm im} \left( A^\top \right) \oplus \ker (A). 
\]
\end{remark}

In our case, we are interested in subspaces \( \langle \mathcal{Q} \rangle \) of \( X(v) \) for which the collection \( \mathcal{Q} \) is feasible. Denote by \( \langle \mathcal{Q} \rangle \) the subspace of \( \Sigma \subseteq \mathbb{R}^N \) spanned by the set of vectors \( \{ \eta^S \mid S \in \mathcal{Q} \} \), and simply write \( \langle S \rangle \) whenever \( \mathcal{Q} = \{S\} \). Let \( A_\mathcal{Q} \) be the subspace of \( X(v) \subseteq \mathbb{R}^N \) defined by \( A_\mathcal{Q} = \{ x \in X(v) \mid x(S) = v(S), \forall S \in \mathcal{Q} \} \).

\medskip

Each projection in \( X(v) \) is naturally decomposed as the sum of two elements: a preimputation and a side payment. Recall that the set of side payments \( \Sigma \) is the linear subspace of \( \mathbb{R}^N \) which is parallel to \( X(v) \) and that for each coalition \( S \in \mathcal{Q} \), the affine subspace \( A_S(v) \) has a corresponding linear subspace \( H_S \) included in \( \Sigma \) with the same dimension as \( A_S(v) \), which was defined by 
\[
H_S \coloneqq \{\sigma \in \Sigma \mid \sigma(S) = 0\}. 
\]
Notice also that \( \eta^S \) is a normal to the subspace \( H_S \). For any subspace \( V \) of a vector space \( E \), the \emph{orthogonal complement} of \( V \), denoted by \( V^\perp \), is defined by 
\[
V^\perp \coloneqq \{x \in E \mid \langle x, y \rangle = 0, \forall y \in V \}. 
\]
Any subspace \( V \) of a finite-dimensional vector space \( E \) induces an orthogonal decomposition of \( E \) into \( E = V \oplus V^\perp \), i.e., for each element \( x \in E \), we can write in a unique way \( x = y + z \) with \( y \in V \) and \( z \in V^\perp \).

\begin{remark}\label{remark: orth-decomp}
Let \( P \) be a projector and \( I \) the identity map. Because the projection of an element is unique, any element \( x \in E \) can be decomposed in a unique way as \( x = y + z \) with \( y = Px \) and \( z = x - Px = (I - P)x\), with \( y \in \text{im}(P) \) and \( z \in {\rm im}(P)^\perp = \ker \left( P^\top \right) \) by Remark~\ref{remark: fundamental-subspaces}. Therefore, we have that \( I - P \) is the projector onto \( \ker \left( P^\top \right) \). 
\end{remark}

Let \( V \) be a linear subspace of \( \Sigma \). Even if it is also a linear subspace of \( \mathbb{R}^N \), in the sequel, we always consider the orthogonal complement \( V^\perp \) as a subset of \( \Sigma \), i.e., 
\[
V^\perp = \{\sigma \in \Sigma \mid \langle \sigma, y \rangle = 0, \, \forall y \in V\}. 
\]
We do not have a well-defined equivalent of the orthogonal complement for affine subspace, but it is still possible to have an orthogonal decomposition. 

\begin{remark} \label{remark: decomposition}
Let \( S \) be a coalition. Using Proposition \ref{prop: first-proj}, any preimputation \( x \) can be decomposed as \( x = \pi_{A_S(v)}(x) - \gamma_S(x) \eta^S \). Naturally, \( \pi_{A_S(v)}(x) \) belongs to \( A_S(v) \) and \( - \gamma_S(x) \eta^S \) belongs to \( \langle S \rangle \), hence we have the decomposition \( X(v) = A_S(v) \oplus \langle S \rangle \). 
\end{remark}

Our main objective in this paper is to compute, for any preimputation outside the core, its projection onto the core. According to the Hilbert Projection Theorem, the projection of \( x \in X(v) \setminus C(v) \) onto the core is the closest preimputation \( y \in C(v) \) from \( x \). Then, \( y \) belongs to a face \( F \) of the core, for which there exists a collection of coalition \( \mathcal{Q} \) such that \( F \subseteq A_\mathcal{Q}(v) \). Then, the first part of our investigation is devoted to the study of the nonemptiness of \( A_\mathcal{Q} \).

\medskip 

In the second part, we construct the projector \( \pi_{\langle \mathcal{Q} \rangle}: \Sigma \to \langle \mathcal{Q} \rangle \), which associates to each side payment \( \sigma \) its projection onto \( \langle \mathcal{Q} \rangle \). The problem amounts to finding the linear combination of the \( \eta^S \), for \( S \in \mathcal{Q} \), which is the closest from \( \sigma \).

\medskip 

From the projector \( \pi_{\langle \mathcal{Q} \rangle} \), we build the projector \( \pi_{A_\mathcal{Q}}: X(v) \to A_\mathcal{Q} \), which associates to each preimputation its projection onto \( A_\mathcal{Q} \), by using the fact that \( \langle \mathcal{Q} \rangle^\perp \) and \( A_\mathcal{Q}  \) are parallel. The projector is expressed in terms of \( \eta^S \) for \( S \in \mathcal{Q} \), with specific coefficients for which an exact formula is provided. 

\medskip 

Finally, we give an alternative formulation of the projector \( \pi_{ A_\mathcal{Q}}: X(v) \to A_\mathcal{Q} \), in terms of the excesses of the given preimputation at the coalitions in \( \mathcal{Q} \), and a set of vectors uniquely defined from the set \( \{ \eta^S \mid S \in \mathcal{Q} \} \).

\subsection{Nonemptiness of \( A_\mathcal{Q} \)} \leavevmode \newline  

In general, the nonemptiness of \( A_\mathcal{Q} \) depends only on \( \mathcal{Q} \).  

\begin{definition}
Let \( \mathcal{Q} \subseteq \mathcal{N} \) be a collection of coalitions. We say that \( \mathcal{Q} \) is an \emph{independent} collection if \( \{ \eta^S \mid S \in \mathcal{Q} \} \) forms a linearly independent set of vectors. 
\end{definition}

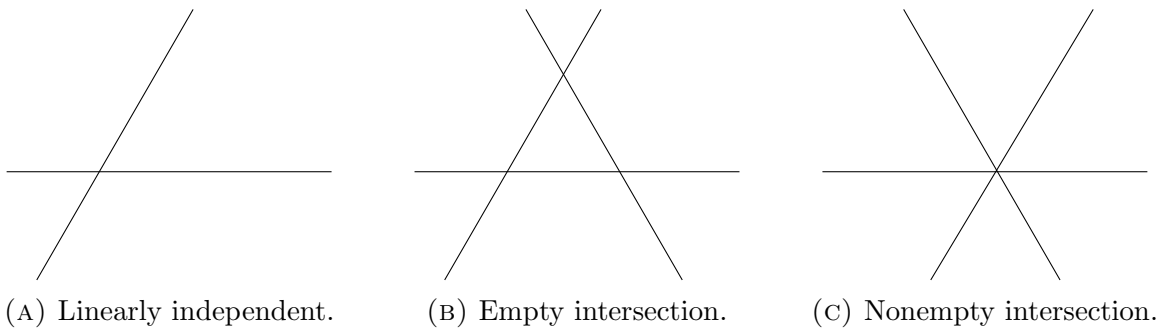
\begin{figure}[ht]
\begin{center}
\begin{subfigure}{0.3\textwidth}
\begin{center}
\begin{tikzpicture}[scale=0.18]
\draw (-7.774, -10) -- (3.774, 10);
\draw (-10, -2) -- (14, -2);
\end{tikzpicture}
\caption{Linearly independent.}
\end{center}
\end{subfigure}
\hspace{0.5cm}
\begin{subfigure}{0.3\textwidth}
\begin{center}
\begin{tikzpicture}[scale=0.18]
\draw (-1.774, 10) -- (9.774, -10);
\draw (-7.774, -10) -- (3.774, 10);
\draw (-10, -2) -- (14, -2);
\end{tikzpicture}
\caption{Empty intersection.}
\end{center}
\end{subfigure}
\hspace{0.5cm}
\begin{subfigure}{0.3\textwidth}
\begin{center}
\begin{tikzpicture}[scale=0.18]
\draw (-6, 10) -- (5.55, -10);
\draw (-4, -10) -- (8, 10);
\draw (-12, -2) -- (12, -2);
\end{tikzpicture}
\caption{Nonempty intersection.}
\end{center}
\end{subfigure}
\caption{Different configurations of hyperplanes in the plane.}
\label{fig: independent}
\end{center}
\end{figure}

A collection which is not independent is said to be \emph{dependent}. We have that \( A_\mathcal{Q} \) is nonempty if \( \mathcal{Q} \) is independent. However, it is not sufficient, as depicted in Figure \ref{fig: independent}. 

\medskip

\begin{definition} 
We say that a coalition function \( v \) is \emph{nonsingular} on a collection \( \mathcal{Q} \) if, for all nonempty dependent subcollections \( \mathcal{T} \subseteq \mathcal{Q} \), we have \( A_\mathcal{T} = \emptyset \). 
\end{definition}

Let \( \mathcal{Q} = \{ S_1, \ldots, S_k \} \subseteq \mathcal{N} \) be a collection of coalitions. We denote by \( H \) and \( G \) respectively the matrices defined by 
\[
H = \begin{bmatrix} \eta^{S_1} & \ldots & \eta^{S_k} \end{bmatrix} \qquad \text{and} \qquad G = H^\top H.
\]
The matrix \( G \) is called the \emph{Gram matrix} of the collection \( \mathcal{Q} \), and its general term satisfies \( g_{ij} = \langle \eta^{S_i}, \eta^{S_j} \rangle \). Then, \( G \) captures all the information about the interdependence and the correlations between the payments of coalitions in \( \mathcal{Q} \). If needed, we can specify the considered collection in the notation by writing \( H_\mathcal{Q} \) and \( G_\mathcal{Q} \). 

\medskip 

The symmetry of the scalar product implies the symmetry of the Gram matrix. Furthermore, for all \( x \in \mathbb{R}^\mathcal{Q} \), we have 
\[
x^\top G x = x^\top H^\top H x = \langle Hx, Hx \rangle = \lVert Hx \rVert^2 \geq 0, 
\]
hence \( G \) is positive semidefinite. We can use the Gram matrix to know whether a collection is independent. 

\begin{proposition}\label{prop: non-singular}
A collection \( \mathcal{Q} \) is independent if and only if \( G \) is nonsingular. 
\end{proposition}

\proof
For all \( x \in \mathbb{R}^\mathcal{Q} \), we have \( x^\top G x \geq 0 \). Then \( G \) is positive semidefinite, i.e., all its eigenvalues are nonnegative. For \( G \) to be nonsingular, we need to have only positive eigenvalues, i.e., to have \( G \) positive definite. Then \( G \) is nonsingular if and only if, for all \( x \in \mathbb{R}^\mathcal{Q} \setminus \{ \vec{o} \, \} \), 
\[
x^\top G x = \lVert Hx \rVert^2 = \left\lVert \sum_{S \in \mathcal{Q}} x_S \eta^S \right\rVert^2 > 0, 
\]
i.e., if and only if, for all \( x \in \mathbb{R}^\mathcal{Q} \setminus \{\vec{o} \, \} \), we have \( \sum_{S \in \mathcal{Q}} x_S \eta^S \neq 0 \), which is the definition of \( \{ \eta^S \mid S \in \mathcal{Q} \} \) being a linearly independent set of vectors. 
\endproof

We now have a sharp characterization of the nonemptiness of \( A_\mathcal{Q} \) when \( v \) is nonsingular on the set system on which it is defined. It is consistent with the fact that our formula for the projector onto \( A_\mathcal{Q} \) involves the inverse of \( G \), as we will see shortly.

\subsection{Construction of the projector \( \pi_{\langle \mathcal{Q} \rangle} : \Sigma \to \langle \mathcal{Q} \rangle \)} \leavevmode \newline 

Because the columns of \( H \) are the normal vectors \( \eta^S \), any element of \( \langle \mathcal{Q} \rangle \) can be written as \( Hy \), with \( y \in \mathbb{R}^\mathcal{Q} \). Intuitively, the side payment between a preimputation \( x \) and its projection onto \( A_\mathcal{Q} \) is a linear combination of the normal vectors \( \eta^S \), for \( S \in \mathcal{Q} \). Now, we are looking for the coefficients of this linear combination. For a given side payment \( \sigma \), we want to find the vector of coefficients \( y \in \mathbb{R}^\mathcal{Q} \) which minimizes
\[
\varphi_H(y) = \lVert \sigma - Hy \rVert^2 = \left( \sigma - Hy \right)^\top \left( \sigma - Hy \right) = \sigma^\top \sigma - y^\top H^\top \sigma - \sigma^\top H y + y^\top G y.
\]
The critical points of \( \varphi_H \) are the elements \( y \in \mathbb{R}^\mathcal{Q} \) satisfying \( \nabla \varphi_H(y) = 0 \). We have the three following gradients: 
\[
\nabla \left( y^\top H^\top \sigma \right) = H^\top \sigma, \qquad \nabla \left( \sigma^\top Hy \right) = H^\top \sigma, \qquad \nabla \left (y^\top G y \right) = 2 G y, 
\]
which lead, by linearity of \( \nabla \), to
\[
\nabla \varphi_H(y) = 2 G y - 2 H^\top \sigma. 
\]
Then, the critical point of \( \varphi_H \) satisfy the so-called \emph{normal equation}, i.e., 
\begin{equation} \label{eq: normal-equation}
Gy = H^\top \sigma. 
\end{equation}
The Hessian of \( \varphi_H \) is \( {\rm Hess}(\varphi_H) = 2G \), which is positive definite when \( \mathcal{Q} \) is independent. Then there is a unique critical point, which is a global minimum. The unique solution \( y \) of the normal equation determines the projection of \( \sigma \) onto \( \langle \mathcal{Q} \rangle \), which is \( Hy \).

\medskip

To manipulate matrix equations such as the normal equations, because not all matrices have an inverse, we need a generalization of usual inverses of matrices.  

\begin{theorem}[\textcite{penrose1955generalized}] \label{th: penrose} \leavevmode \newline
For any matrix \( A \), there exists a unique matrix \( X \) satisfying
\[
AXA = A, \qquad XAX = X, \qquad \left( AX \right)^\top = AX, \qquad \text{and} \qquad \left( XA \right)^\top = XA. 
\]
\end{theorem}

These equations are called the \emph{Penrose equations}, and their solution is called the \emph{Moore-Penrose inverse}, or \emph{generalized inverse}, of \( A \) and is denoted by \( A^\dag \). We notice that if \( A \) is nonsingular, we have that \( A^\dag = A^{-1} \). Moreover, if the columns of \( A \) are linearly independent, we have 
\[
A^\dag = \left( A^\top A \right)^{-1} A^\top, 
\]
and if the rows of \( A \) are linearly independent, we have
\[
A^\dag = A^\top \left( A A^\top \right)^{-1}. 
\]
When the columns of \( A \) are linearly independent, the matrix \( A^\dag \) is a \emph{left inverse} for \( A \), because \( A^\dag A = \left( A^\top A \right)^{-1} A^\top A = {\rm Id} \). Similarly, when the rows of \( A \) are linearly independent, \( A^\dag \) is a \emph{right inverse}. Then, if \( \mathcal{Q} \) is independent, the Moore-Penrose inverse of \( H \) is the matrix 
\[
H^\dag = \left( H^\top H \right)^{-1} H^\top = G^{-1} H^\top. 
\]
Then, we can solve the normal equation \eqref{eq: normal-equation}, and find \( \pi_{\langle \mathcal{Q} \rangle}(\sigma) \):
\[
y = G^{-1} H^\top \sigma = H^\dag \sigma \qquad \text{and} \qquad \pi_{\langle \mathcal{Q} \rangle}(\sigma) = H y = H H^\dag \sigma. 
\]

\subsection{Construction of the projector $\pi_{A_\mathcal{Q}}: X(v) \to A_\mathcal{Q}$} \leavevmode \newline

We have built the projector onto \( {\rm im}(H) = \langle \mathcal{Q} \rangle \), from the matrix \( H \) and its Moore-Penrose inverse \( H^\dag \). The initial goal is to build a projector onto the affine subspace \( A_\mathcal{Q} \) for all preimputations. We introduce \( b^\mathcal{Q} \in \mathbb{R}^\mathcal{Q} \) defined, for all \( S \in \mathcal{Q} \), by 
\[
b^\mathcal{Q}_S = v(S) - \frac{\lvert S \rvert}{\lvert N \rvert} v(N). 
\]
Let \( L = H^\top \). We can describe \( A_\mathcal{Q} \) using the matrix \( L \) and \( b^\mathcal{Q} \), by
\[
A_\mathcal{Q} = \left\{ x \in X(v) \mid Lx = b^\mathcal{Q} \right \}. 
\]
The linear subspace of \( \Sigma \) parallel to \( A_\mathcal{Q} \) is therefore \( \{ \sigma \in \Sigma \mid L \sigma = 0 \} = \ker (L) \cap \Sigma \). From Remark~\ref{remark: fundamental-subspaces}, we know that \( \mathbb{R}^N = \langle \mathcal{Q} \rangle \oplus \ker (L) \). Then, \( \pi_{\ker (L)} = {\rm Id} - \pi_{\langle \mathcal{Q} \rangle} \). To deduce \( \pi_{A_\mathcal{Q}} \) from \( \pi_{\ker (L)} \), we use the following lemma. 

\begin{lemma}[\textcite{bauschke2011convex}]\label{lemma: proj-affine} \leavevmode \newline
Let \( K \) be a nonempty closed convex subset of \( \mathbb{R}^N \), and let \( x, y \in \mathbb{R}^N \). Then
\[
\pi_{y + K}(x) = y + \pi_K (x-y). 
\]
\end{lemma}

Using Lemma \ref{lemma: proj-affine}, we have, for any \( y \in A_\mathcal{Q} \), that 
\[
\pi_{A_\mathcal{Q}}(x) = y + \pi_{\ker(L)}(x - y) = y + \left({\rm Id} - H H^\dag \right) \left( x - y \right) = x - H G^{-1} L \left( x - y \right). 
\]
Indeed, any element of \( A_\mathcal{Q} \) can be written as \( y + z \) with \( z \in \ker (L) \) because any side payment in \( \ker (L) \) preserves the payment of \( S \in \mathcal{Q} \). Because \( y \in A_\mathcal{Q} \), we have \( Ly = b^\mathcal{Q} \). Denote by \( e_\mathcal{Q} (x) \in \mathbb{R}^\mathcal{Q} \) the vector defined, for all \( S \in \mathcal{Q} \), by \( \left[e_\mathcal{Q}(x)\right]_S = e_S(x) \).

\medskip 

We are now able to state the main result of this paper. 

\begin{theorem}\label{th: main-formula-proj}
Let \( \mathcal{Q} \) be an independent collection. For all \( x \in X(v) \), we have
\[
\pi_{A_\mathcal{Q}}(x) = x + \sum_{S \in \mathcal{Q}} \gamma^\mathcal{Q}_S(x) \eta^S, \qquad \text{where} \quad \gamma^\mathcal{Q}_S(x) = \left[ G^{-1} e_\mathcal{Q}(x) \right]_S.
\]
 
\end{theorem}

\proof
To prove this theorem, we use the Projection Theorem for affine subspaces. Let $x$ be a preimputation. First, let us show that \( \pi_{A_\mathcal{Q}}(x) \in X(v) \). We have that \( \pi_{A_\mathcal{Q}}(x)(N) = x(N) + \sum_{S \in \mathcal{Q}} \gamma^\mathcal{Q}_S(x) \eta^S(N) = v(N) \). Also, notice that 
\[
L(x - y) = Lx - b^\mathcal{Q} = - e_\mathcal{Q}(x). 
\]
Now, let us prove that \( \pi_{A_\mathcal{Q}}(x) \) belongs to \( A_\mathcal{Q} \), i.e., that \( L \pi_{A_\mathcal{Q}}(x) = b^\mathcal{Q} \):
\[
L \pi_{A_\mathcal{Q}}(x) = L \left( x + H G^{-1} e_\mathcal{Q}(x) \right) = Lx - L H G^{-1} \left( Lx - b^\mathcal{Q} \right). 
\]
Because \( L= H^T \), we have that \( LH = G \), and it follows that
\[
L \pi_{A_\mathcal{Q}}(x) = Lx - \left( Lx - b^\mathcal{Q} \right) = b^\mathcal{Q}, 
\]
therefore \( \pi_{A_\mathcal{Q}}(x) \in A_\mathcal{Q} \). Let us prove now that \( x - \pi_{A_\mathcal{Q}}(x) \) is orthogonal to \( A_\mathcal{Q} \). Denote by \( \gamma^\mathcal{Q}_S(x) \) the coordinate \( \left[ G^{-1} e_\mathcal{Q}(x) \right]_S \). For all \( y \in A_\mathcal{Q} \), we have 
\[
- \langle y - \pi_{A_\mathcal{Q}}(x), x - \pi_{A_\mathcal{Q}}(x) \rangle = \sum_{S \in \mathcal{Q}} \gamma^\mathcal{Q}_S(x) \langle y - \pi_{A_\mathcal{Q}}(x), \eta^S \rangle.
\]
Because both \( y \) and \( \pi_{A_\mathcal{Q}}(x) \) belong to \( A_\mathcal{Q} \), we have that \( y(N) = v(N) = \pi_{A_\mathcal{Q}}(x)(N)\) and, for all \( S \in \mathcal{Q} \), that \( y(S) = \pi_{A_\mathcal{Q}}(x)(S) \). Then, \( \langle y, \eta^S \rangle = \langle \pi_{A_\mathcal{Q}}(x), \eta^S \rangle \), and finally 
\[
- \langle y - \pi_{A_\mathcal{Q}}(x), x - \pi_{A_\mathcal{Q}}(x) \rangle = \sum_{S \in \mathcal{Q}} \gamma^\mathcal{Q}_S(x) \left( \langle y, \eta^S \rangle - \langle \pi_{A_\mathcal{Q}}(x), \eta^S \rangle \right) = 0, 
\]
which concludes the proof by the Projection Theorem. 
\endproof

We have now found a closed-form formula for the projection of any given preimputation \( x \) onto an affine subspace \( A_\mathcal{Q} \), provided that \( \mathcal{Q} \) is an independent collection. It requires no iterative computations, and can be implemented within a few lines of codes. The main algorithmic effort is to compute the inverse the Gram matrix \( G \). One way to avoid that is to use an existing solver of linear systems to find \( \gamma^\mathcal{Q}(x) \) satisfying 
\begin{equation} \label{eq: cramer}
G \gamma^\mathcal{Q}(x) = e_\mathcal{Q}(x). 
\end{equation}

\begin{remark}[Cramer's rule \parencite{cramer1750introduction}] \label{remark: cramer-rule}
We can use Cramer's rule to solve the linear system~\eqref{eq: cramer}, which express the coordinates of \( \gamma^S(x) \) in terms of the determinant of \( G \), called the \emph{Gramian} of \( \mathcal{Q} \). Let \( G^S_x \) denote the matrix form where we replace the column composed of all the scalar products involving \( \eta^S \) with the vector \( e_\mathcal{Q}(x) \), i.e., 
\[
G^{S_i}_x = \begin{pmatrix} \lVert \eta^{S_1} \rVert^2 & \ldots & \langle \eta^{S_1}, \eta^{S_{i-1}} \rangle & e(S_1, x) & \langle \eta^{S_1}, \eta^{S_{i+1}} \rangle & \ldots & \langle \eta^{S_1}, \eta^{S_k} \rangle \\
\langle \eta^{S_2}, \eta^{S_1} \rangle & \ldots & \langle \eta^{S_2}, \eta^{S_{i-1}} \rangle & e(S_2, x) & \langle \eta^{S_2}, \eta^{S_{i+1}} \rangle & \ldots & \langle \eta^{S_2}, \eta^{S_k} \rangle \\
\vdots & & & \vdots & & & \vdots \\
\langle \eta^{S_j}, \eta^{S_1} \rangle & \ldots & \langle \eta^{S_j}, \eta^{S_{i-1}} \rangle & e(S_j, x) & \langle \eta^{S_j}, \eta^{S_{i+1}} \rangle & \ldots & \langle \eta^{S_j}, \eta^{S_k} \rangle \\
\vdots & & & \vdots & & & \vdots \\
\langle \eta^{S_k}, \eta^{S_1} \rangle & \ldots & \langle \eta^{S_k}, \eta^{S_{i-1}} \rangle & e(S_k, x) & \langle \eta^{S_k}, \eta^{S_{i+1}} \rangle & \ldots & \lVert \eta^{S_k} \rVert^2 \end{pmatrix} 
\]
Cramer's rule gives that 
\[
\gamma^\mathcal{Q}_S(x) = \frac{\det G^S_x}{\det G}, 
\]
therefore the formula of the projector for any preimputation \( x \) is 
\[
\pi_{A_\mathcal{Q}}(x) = x + \frac{1}{\det G} \sum_{S \in \mathcal{Q}} \det G^S_x \cdot \eta^S. 
\]
The relevance of Cramer's rule here is mainly for theoretical considerations. Even if the computational complexity of Cramer's rule has greatly reduced in the past years (see \textcite{habgood2012condensation}), the Cholesky decomposition is preferred to solve linear systems of equations with a symmetric positive definite matrix of coefficients. 
\end{remark}

\begin{example}
Let \( S \) be a coalition. We apply Theorem~\ref{th: main-formula-proj} and Cramer's rule to compute the projector of \( x \in X(v) \) onto \( A_S \). We have that \( G = \begin{pmatrix} \lVert \eta^S \rVert^2 \end{pmatrix} \), and \( G^S_x = \begin{pmatrix} e_S(x) \end{pmatrix} \). It follows that 
\[
\pi_{A_S}(x) = x + \frac{1}{\det G} \sum_{S \in \{S\}} \det G^S_x \cdot \eta^S = x + \frac{e_S(x)}{\lVert \eta^S \rVert^2} \eta^S,
\]
which coincides with the formula of Proposition \ref{prop: first-proj}.
\end{example}

\subsection{Alternative formula for \( \pi_{A_\mathcal{Q}}: X(v) \to A_\mathcal{Q} \)} \leavevmode \newline 

It is possible to express the projection of a given preimputation in a simpler way, using a new set of vectors, uniquely determined from \( \{ \eta^S \mid S \in \mathcal{Q} \} \). Denote by \( \delta_{ij} \) the Kronecker delta, defined by \( \delta_{ij} = 1 \) if \( i = j \) and \( \delta_{ij} = 0 \) otherwise. 

\begin{definition}
Let \( \{ e_1, \ldots, e_k \} \) and \( \{ f_1, \ldots, f_k \} \) be two sets of vectors in \( \mathbb{R}^N \). We say that they form a \emph{biorthogonal system} if they satisfy \( \langle e_i, f_j \rangle = \delta_{ij} \). 
\end{definition}

A common example of a biorthogonal system is the pair formed by a basis of a Euclidean space together with its dual basis. Then, for simplicity, if the sets are linearly independent, we say that \( \{f_1, \ldots, f_k\} \) is a dual basis of \( \{e_1, \ldots, e_k\} \), as it is the case in \( {\rm span}(\{e_1, \ldots, e_k\}) \). 

\begin{proposition}\label{prop: dual-basis}
Let \( \{e_1, \ldots, e_k\} \) be a set of linearly independent vectors in \( \mathbb{R}^N \), and let \( \{ f_1, \ldots, f_k \} \) be a dual basis of \( \{e_1, \ldots, e_k\} \). Define two \((\lvert N \rvert \times k)\)-matrices \( E \) and \( F \) by \( E = \begin{pmatrix} e_1 & \ldots & e_k \end{pmatrix} \) and \( F = \begin{pmatrix} f_1 & \ldots & f_k \end{pmatrix} \), as well as \( G^E = E^\top E \) and \( G^F = F^\top F \). The following equalities hold. 
\[
E^\top F = F^\top E = {\rm Id}_k, \qquad \text{and} \qquad \left( G^E \right)^{-1} = G^F.
\]
\end{proposition}

\proof
Let us look at the entries of \( E^\top F \). From the definition of biorthogonal systems, we deduce \( \left( E^\top F \right)_{ij} = \langle e_i, f_j \rangle = \delta_{ij} \), and \( \left( F^\top E \right)_{ij} = \langle f_i, e_j \rangle = \delta_{ij} \), then \( E^\top F = {\rm Id}_k = F^\top E \). Let us prove now that \( \left( G^E \right)^{-1} = G^F \). Denote \( K = \{1, \ldots, k\} \) and let \( x \in {\rm span}(\{e_1, \ldots, e_k\}) \). We write the coordinates of \( x \) in the bases \( \{e_1, \ldots, e_k\} \) and \( \{f_1, \ldots, f_k\} \) respectively by \( x = \sum_{i \in K} \alpha_i e_i \) and \( x = \sum_{i \in K} \beta_i f_i \). We denote by \( g^E_{ij} \) and \( g^F_{ij} \) the general terms of \( G^E \) and \( G^F \) respectively. We can express the coordinates in one basis in terms of the coordinates in the other basis by 
\[ \begin{aligned}
\alpha_p & = \sum_{i \in K} \alpha_i \langle e_i, f_p \rangle = \langle x, f_p \rangle = \sum_{i \in K} \beta_i \langle f_i, f_p \rangle = \sum_{i \in K} g^F_{ip} \beta_i, \qquad \forall p \in K, \\
\text{and} \qquad \beta_q & = \sum_{j \in K} \beta_j \langle f_j, e_q \rangle = \langle x, e_q \rangle = \sum_{j \in K} \alpha_j \langle e_j, e_q \rangle = \sum_{j \in K} g^E_{jq} \alpha_j, \qquad \forall q \in K.
\end{aligned} \]
Combining these two expressions yields 
\[
\alpha_p = \sum_{i \in K} g^F_{ip} \left( \sum_{j \in K} g^E_{ji} \alpha_j \right) = \sum_{j \in K} \left( \sum_{i \in K} g^E_{ji} g^F_{ip} \right) \alpha_j = \sum_{j \in K} \left( G^E G^F \right)_{jp} \alpha_j, \qquad \forall p \in K. 
\]
Then, we have that \( \left( G^E G^F \right)_{ij} = \delta_{ij} \), and doing the same calculations on the \( \beta \)'s gives \( \left( G^F G^E \right)_{ij} = \delta_{ij} \), hence \( G^F \) is the inverse of \( G^E \). 
\endproof

Let \( \mathcal{Q} = \{ S_1, \ldots, S_k \} \) be an independent collection. We denote by \( \{ h^{S_1}, \ldots, h^{S_k} \} \) the dual basis of \( \{ \eta^{S_1}, \ldots, \eta^{S_k} \} \). Set \( H^\circ = \begin{pmatrix} h^{S_1} & \ldots & h^{S_k} \end{pmatrix} \). 

\begin{lemma} \label{lemma: dual-basis}
Let \( \mathcal{Q} \) be an independent collection of coalitions. Then 
\[
h^{S_j} = \sum_{i \in K} g^{[-1]}_{ij} \eta^{S_i} \qquad \text{and} \qquad H^\circ = H G^{-1} = L^\dag, 
\]
where \( g^{[-1]}_{ij} \) denotes the general term of \( G^{-1} \). 
\end{lemma}

\proof
The second fact is simply the translation of the first one in terms of matrices, we will therefore only prove the first fact. Let \( x \in \langle \mathcal{Q} \rangle \), and we write the coordinates of \( x \) in basis \( \{ \eta^S \mid S \in \mathcal{Q} \} \) and \( \{ h^S \mid S \in \mathcal{Q} \} \) respectively by \( x = \sum_{S \in \mathcal{Q}} \alpha_S \eta^S \) and \( x = \sum_{S \in \mathcal{Q}} \beta_S h^S \). Using Proposition \ref{prop: dual-basis} gives
\[
\alpha_{S_p} = \sum_{i \in K} \alpha_{S_i} \langle \eta^{S_i}, h^{S_p} \rangle = \langle x, h^{S_p} \rangle = \sum_{i \in K} \beta_{S_i} \langle h^{S_i}, h^{S_p} \rangle = \sum_{i \in K} g^{[-1]}_{ip} \beta_{S_i}, \qquad \forall p \in K.
\]
By setting all the coefficients \( \beta_S \) to \( 0 \) except for \( \beta_{S_j} = 1 \), we get \( \alpha_{S_p} = g^{[-1]}_{jp} \) for all \( p \in K \). Then, \( h^{S_j} = x = \sum_{i \in K} g^{[-1]}_{ij} \eta^{S_i} \). We can rewrite it as 
\[
h^{S_j} = \eta^{S_1} g^{[-1]}_{1j} + \eta^{S_2} g^{[-1]}_{2j} + \ldots + \eta^{S_k} g^{[-1]}_{kj} = \begin{bmatrix} \eta^{S_1} & \ldots & \eta^{S_k} \end{bmatrix} \left( G^{-1} \right)^{\text{col}}_j,
\]
with \( \left( G^{-1} \right)^{\text{col}}_j \) denoting the \( j \)-th column of \( G^{-1} \), and then we have \( H^\circ = H G^{-1} \). Because the rows of \( L \) are linearly independent, we have 
\[
L^\dag = L^\top \left( L L^\top \right)^{-1} = H \left( H^\top H \right)^{-1} = H G^{-1} = H^\circ, 
\]
which concludes the proof. 
\endproof

Using Lemma~\ref{lemma: dual-basis}, we can rewrite the formula of the projector. 

\begin{theorem}\label{th: proj-h}
Let \( \mathcal{Q} \) be an independent collection. For all \( x \in X(v) \), we have 
\[
\pi_{A_\mathcal{Q}}(x) = x + \sum_{S \in \mathcal{Q}} e_S(x) h^S.
\]
\end{theorem}

\proof
From Theorem~\ref{th: main-formula-proj}, we already have that \( \pi_{A_\mathcal{Q}}(x) = x + H G^{-1} e_\mathcal{Q}(x) \). By Lemma~\ref{lemma: dual-basis}, we have \( H G^{-1} = H^\circ \), and therefore, 
\[
\pi_{A_\mathcal{Q}}(x) = x + H^\circ e_\mathcal{Q}(x) = x + \sum_{S \in \mathcal{Q}} e_S(x) h^S. 
\]
\endproof

\begin{figure}[ht]
\begin{center}
\begin{subfigure}{0.49\textwidth}
\begin{center}
\begin{tikzpicture}[scale=0.2]
\draw[->] (-8, 6) -- (-8, 8) node[above]{\( \eta^T \)};
\draw[->] (2.25, 9.034) -- (4.25, 10.184) node[above right]{\( \eta^S \)};

\draw (1.774, 10) node[above left]{\( A_S \)} -- (11.774, -10);
\draw (-10, 6) node[left]{\( A_T \)} -- (14, 6);

\filldraw[black] (-5,-5) circle (5pt) node[below left]{\( x \)};

\draw[dashed] (-5, -5) -- (3.774, 6);

\draw[thick, blue, ->] (-5, -5) -- node [midway, below right]{\( \gamma^\mathcal{Q}_S(x) \eta^S \)} (3.774, -0.613);

\draw[color=white!50, fill=white!50] (-2.3, 1.5) rectangle (3.5, 3.8);
\draw[thick, blue, ->] (3.774, -0.613) -- node[midway, left]{\( \gamma^\mathcal{Q}_T(x) \eta^T \)} (3.774, 6);
\end{tikzpicture}
\end{center}
\end{subfigure}
\begin{subfigure}{0.49\textwidth}
\begin{center}
\begin{tikzpicture}[scale=0.2]
\draw[->] (-8, 6) -- (-8, 8) node[above]{\( \eta^T \)};
\draw[->] (2.25, 9.034) -- (4.25, 10.184) node[above right]{\( \eta^S \)};

\draw (1.774, 10) node[above left]{\( A_S \)} -- (11.774, -10);
\draw (-10, 6) node[left]{\( A_T \)} -- (14, 6);

\filldraw[black] (-5,-5) circle (5pt) node[below left]{\( x \)};

\draw[dashed] (-5, -5) -- (3.774, 6);

\draw[thick, blue, ->] (-5, -5) -- node [midway, below]{\footnotesize \( e_T (x) h^T \)} (9.274, -5);

\draw[thick, blue, ->] (9.274, -5) -- node[midway, above right]{\footnotesize \( e_S(x) h^S \)} (3.774, 6);
\end{tikzpicture}
\end{center}
\end{subfigure}
\caption{Decompositions of the side payment between \( x \) and its projection onto \( A_\mathcal{Q} \) with \( \mathcal{Q} = \{S, T\} \). The vectors \( \eta^S\), \( \eta^T \), \( h^S \) and \( h^T \) are not properly scaled.}
\label{fig: two-proj}
\end{center}
\end{figure}
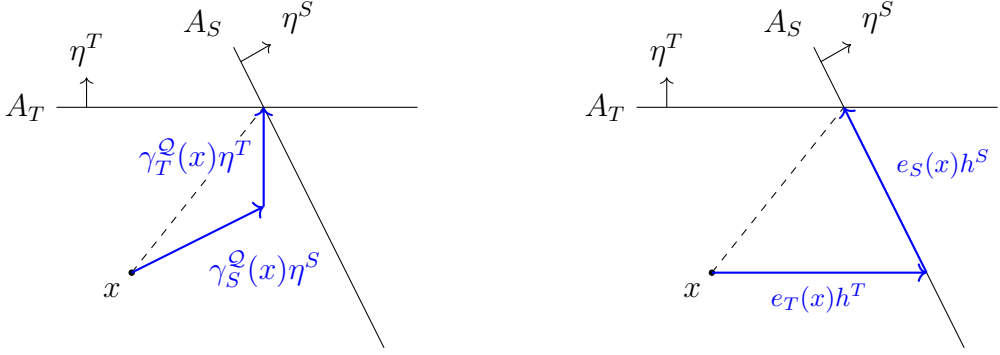

\medskip

It is possible to obtain a simple matrix as a projector using an orthonormal basis of \( \langle \mathcal{Q} \rangle \). To do so, we need to be able to project a vector \( \eta^S \) onto the line \( \langle T \rangle \) spanned by another vector \( \eta^T \), with \( S \) and \( T \) being two coalitions. 

\begin{lemma} \label{lemma: proj-gram-schmidt}
Let \( S \) be a coalition, and let \( \sigma \in \Sigma \) be a side payment. Then the projection of \( \sigma \) onto \( \langle S \rangle \) is 
\[
\pi_{\langle S \rangle}(\sigma) = \frac{\sigma(S)}{\lVert \eta^S \rVert^2} \eta^S. 
\]
\end{lemma}

\proof
Let \( (N, o) \) be the null game, assigning to every coalition a worth of \( 0 \). We have that \( \Sigma = X(o) \) and \( A_S(o) = H_S \). Then, the map \( \pi_{A_S(o)} \) projects a side payment onto \( H_S \). Using Remark~\ref{remark: fundamental-subspaces}, we know that 
\[
\sigma = \pi_{H_S}(\sigma) + \pi_{\langle S \rangle}(\sigma).
\] 
Using Proposition~\ref{prop: first-proj}, we have that 
\[
\pi_{H_S}(\sigma) = \pi_{A_S(o)}(\sigma) = \sigma + \frac{o(S) - \sigma(S)}{\lVert \eta^S \rVert^2} \eta^S = \sigma - \frac{\sigma(S)}{\lVert \eta^S \rVert^2} \eta^S, 
\]
from which we deduce that \( \pi_{\langle S \rangle}(\sigma) = \frac{\sigma(S)}{\lVert \eta^S \rVert^2} \eta^S \). 
\endproof

\begin{remark}
Let \( \sigma \) be a side payment and \( S \) be a coalition. Then, we have that \( \sigma(N) = 0 \), which leads to \( \langle \sigma, \eta^S \rangle = \langle \sigma, \mathbf{1}^S \rangle - \frac{\lvert S \rvert}{\lvert N \rvert} \langle \sigma, \mathbf{1}^N \rangle = \langle \sigma, \mathbf{1}^S \rangle = \sigma(S) \). We can therefore retrieve the well-known formula: 
\[
\pi_{\langle S \rangle}(\sigma) = \frac{\langle \sigma, \eta^S \rangle}{\lVert \eta^S \rVert^2} \eta^S. 
\]
\end{remark}

We can therefore use the modified Gram-Schmidt process to obtain an orthonormal basis of \( \langle \mathcal{Q} \rangle \), using Lemma~\ref{lemma: proj-gram-schmidt}. Let \( \mathcal{V} \) denote the orthonormal basis of \( \langle \mathcal{Q} \rangle \) obtained via Gram-Schmidt process. Denote \( V \coloneqq \begin{pmatrix} v_1 & \ldots & v_k \end{pmatrix} \). The matrix \( V \) is therefore semi-orthogonal, i.e., \( V^\dag = V^\top \) and \( V^\top V = {\rm Id} \). The projector onto \( \langle \mathcal{Q} \rangle \) can be reformulated as \( \pi_{\langle \mathcal{Q} \rangle} = V V^\top \). Then, the projection of a preimputation \( x \) onto \( A_\mathcal{Q} \) is, for all \( y \in A_\mathcal{Q} \), 
\[
\pi_{A_\mathcal{Q}}(x) = x - V V^\top (x-y).
\]
Such an element \( y \in A_\mathcal{Q} \) can be found using the Moore-Penrose inverse of \( L \).

\section{Using the appropriate projector on a given preimputation}\label{sec: algo-proj}

In this section, we finalize the construction of the algorithm computing the projection of a given preimputation onto the core. This far, we know how to project onto the intersection of severable affine subspaces \( A_S(v) \), as long as the coalitions form an independent collection. The question now is to identify, for a given preimputation, the independent collection that we use to build its projector. 

\medskip 

Assume now that \( \mathcal{Q} \) is feasible, but not independent. Let \( (N, v) \) be a balanced game that is nonsingular on \( \mathcal{Q} \). Then, the affine subspace \( A_\mathcal{Q} \) is empty, and there is no obvious choice of independent subcollection \( \mathcal{T} \subseteq \mathcal{Q} \) such that, for all \( x \in X_\mathcal{Q} \), we have \( \pi_{A_\mathcal{T}}(x) \in C(v) \). Moreover, there may be several independent subcollections satisfying this property, and they might be nonmaximal with respect to inclusion. 

\medskip

The next result is the first step to the construction of an algorithm finding a proper independent subcollection \( \mathcal{T} \subseteq \mathcal{Q} \), and generalize Proposition~\ref{prop: chi} to the general case. Let \( \mathcal{Q} \) be an independent collection, \( T \not \in \mathcal{Q} \) be a coalition and \( x \in X(v) \). We define
\[
\chi_\mathcal{Q}(T, x) = \sum_{S \in \mathcal{Q}} \left[ \langle \eta^S, \eta^T \rangle \cdot \det G^S_x - e_S(x) \cdot \det G \right].
\]

\begin{proposition}
Let \( \mathcal{Q} \subseteq \mathcal{N} \) be independent, \( T \in \mathcal{N} \setminus \mathcal{Q} \) and \( x \in X(v) \). Then 
\[
\pi_{A_\mathcal{Q}}(x) \in A^\geq_T \quad \text{ if and only if } \quad \chi_\mathcal{Q}(T, x) \geq 0. 
\]
\end{proposition}

\proof
First, we study the excess of \( T \) at the projection onto \( A_\mathcal{Q} \):
\[ 
e_T(\pi_{A_\mathcal{Q}}(x)) = v(T) - x(T) - \sum_{S \in \mathcal{Q}} \gamma^\mathcal{Q}_S(x) \eta^S(T).
\]
Using Cramer's rule (see Remark~\ref{remark: cramer-rule}) gives
\[
e_T(\pi_{A_\mathcal{Q}}(x)) = e_T(x) - \sum_{S \in \mathcal{Q}} \frac{\det G^S_x}{\det G} \langle \eta^S, \eta^T \rangle = -\left( \det G \right)^{-1} \chi_\mathcal{Q}(T, x). 
\]
The projection lies into \( A^\geq_T \) if and only if \( e_T(\pi_{A_\mathcal{Q}}(x)) \) is nonpositive, and therefore if and only if \( \chi_\mathcal{Q}(T, x) \) is nonnegative. 
\endproof

\begin{figure}[ht]
\begin{center}
\begin{tikzpicture}[scale=0.3]
\fill[gray!15] (5.1548, -2) -- (1, 5.1957) -- (3.774, 10) -- (11.3677, 10) -- (14, 5.4411) -- (14, -2) -- cycle;

\draw (-1.774, 10) node[above left]{\( A_S \)} -- (9.774, -10);
\draw (-7.774, -10) -- (3.774, 10) node[above right]{\( A_T \)};
\draw (-10, -2) -- (14, -2) node[right]{\( A_K \)};

\path (7, 4) node[below] {\( C(v) \)};

\filldraw[black] (-6,-5) circle (5pt) node[below left]{\( x \)};

\draw[dashed] (-6, -5) -- (3.665, 0.5804);

\draw[thick, blue, ->] (-6, -5) -- (1, 5.196);
\draw[thick, blue, ->] (-6, -5) -- (-3.155, -2);
\draw[thick, blue, ->] (-6, -5) -- (5.155, -2);
\end{tikzpicture}
\caption{Projections on \( A_\mathcal{T} \) for all maximal independent collection \( \mathcal{T} \subseteq \mathcal{Q} = \{S, T, K \} \), with the projection (dashed) on the core.}
\label{fig: multi-proj}
\end{center}
\end{figure}
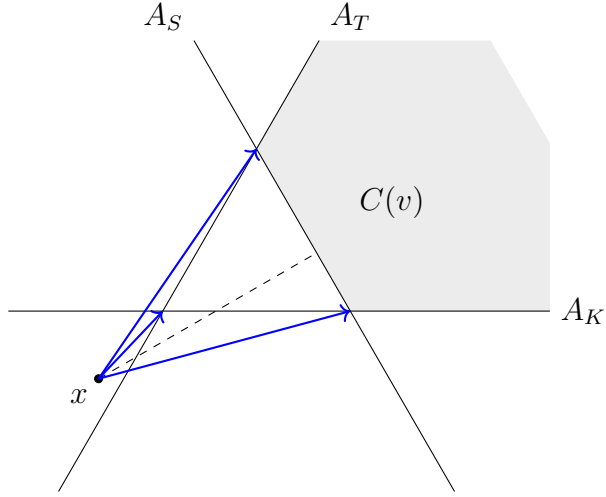

We now have a characterization to check whether the projection onto an affine subspace \( A_\mathcal{T} \) with \( \mathcal{T} \subseteq \mathcal{Q} \) cross or lies into all the affine hyperplanes \( A_S \), for \( S \in \mathcal{Q} \). 

\medskip 

One drawback of this characterization is the frequent use of determinants. But the computation of \( \det G^S_x \) can be done in terms of Gramians using Laplace's expansion formula. To compute determinant of these Gramians, we use the following result. 

\begin{proposition}\label{prop: gramian}
Let \( \mathcal{Q} \) be an independent collection, and \( \mathcal{V} = \{v_S \mid S \in \mathcal{Q} \} \) be an orthonormal basis of \( \langle \mathcal{Q} \rangle \). Then, 
\[
\det G = \left( \prod_{S \in \mathcal{Q}} v_S(S) \right)^2. 
\]
\end{proposition}

\proof
The proof is based on the QR decomposition of \( H \) by the Gram-Schmidt process. Indeed, the QR decomposition of \( H = QR \) gives an orthogonal matrix \( Q \) whose columns form an orthonormal basis \( \mathcal{V} = \{v_{S_1}, \ldots, v_{S_k} \} \) of \( \langle \mathcal{Q} \rangle \), and an upper triangular matrix \( R \) defined by 
\[
R = \begin{pmatrix} \left\langle v_1, \eta^{S_1} \right\rangle & \left\langle v_1, \eta^{S_2} \right\rangle & \left\langle v_1, \eta^{S_3} \right\rangle & \ldots & \left\langle v_1, \eta^{S_k} \right\rangle \\
0 & \left\langle v_2, \eta^{S_2} \right\rangle & \left\langle v_2, \eta^{S_3} \right\rangle & \ldots & \left\langle v_2, \eta^{S_k} \right\rangle \\
0 & 0 & \left\langle v_3, \eta^{S_3} \right\rangle & \ldots & \left\langle v_3, \eta^{S_k} \right\rangle \\
\vdots & \vdots & \vdots & \ddots & \vdots \\
0 & 0 & 0 & \ldots & \left\langle v_k, \eta^{S_k} \right\rangle \end{pmatrix}.
\]
Therefore, the determinant of \( R \) is \( \det R = \prod_{S \in \mathcal{Q}} \langle v_S, \eta^S \rangle = \prod_{S \in \mathcal{Q}} v_S(S) \) because, for all \( S \in \mathcal{Q} \), \( v_S \in \langle \mathcal{Q} \rangle \subseteq \Sigma \). By the multiplicativity of the  determinant, we have that 
\[
\det G = \det \left( H^\top H \right) = \det \left( R^\top Q^\top Q R \right) = \left( \det R \right)^2, 
\]
which leads to \( \det G = \left( \prod_{S \in \mathcal{Q}} v_S(S) \right)^2 \). 
\endproof

\begin{algorithm}
\caption{Gramian \( \det G \) computation} \label{algo: gramian}
\begin{algorithmic}[1]
\Require The Gramian \( \det G \) of \( \mathcal{Q} \), an orthonormal basis \( \mathcal{V} \) of \( \langle \mathcal{Q} \rangle \), a coalition \( S \not \in \mathcal{Q} \)
\Ensure The updated Gramian of \( \mathcal{Q} \cup \{S\} \), an updated orthonormal basis of \( \mathcal{Q} \cup \{S\} \)
\Procedure{UpdateGramian}{$\mathcal{Q}, S$}
\State Set \( v_S \gets \eta^S \)
\For{\( v \in \mathcal{V} \)}
\State \( v_S \gets v_S - v(S) \cdot v \) \Comment{see Lemma~\ref{lemma: proj-gram-schmidt}}
\EndFor
\If{\( \lVert v_S \rVert \neq 0 \)}
\State \( v_S \gets v_S / \lVert v_S \rVert \)
\State \( \det G \gets \left(v_S (S)\right)^2 \cdot \det G \) \Comment{see Proposition~\ref{prop: gramian}}
\State \( \mathcal{V} \gets \mathcal{V} \cup \{v_S \} \)
\EndIf
\State \textbf{return} \( \mathcal{V} \), \( \det G \)
\EndProcedure
\end{algorithmic}
\end{algorithm}

Combining this result with the Gram-Schmidt process leads to Algorithm~\ref{algo: gramian}, which has a complexity of \( \mathcal{O}(n k^2) \), which is slightly slower than the usual ones used to compute determinants, which have a complexity of \( \mathcal{O}(k^3) \). However, the idea is not to compute directly a Gramian with this algorithm, but to derive the Gramian of a collection \( \mathcal{Q} \) by updating the one of the collection \( \mathcal{Q} \setminus S \) with \( S \in \mathcal{Q} \). In the next result, we present a projector onto lines in \( \Sigma \) that we used in the algorithm computing the Gramian. 

\begin{lemma}\label{lemma: proj-gram-schmidt}
Let \( S \) be a coalition, and let \( \sigma \in \Sigma \) be a side payment satisfying \( \lVert \sigma \rVert = 1 \). Then the projection of \( \eta^S \) onto \( {\rm span}(\sigma) \) is 
\[
\pi_{{\rm span}(\sigma)}(\eta^S) = \sigma(S) \cdot \sigma. 
\]
\end{lemma}

\proof
The usual projector onto lines gives
\[
\pi_{{\rm span}(\sigma)}(\eta^S) = \frac{\langle \eta^S, \sigma \rangle}{\langle \sigma, \sigma \rangle} \sigma = \langle \eta^S, \sigma \rangle \sigma. 
\]
Because \( \sigma \) is a side payment, we have \( \langle \mathbf{1}^N, \sigma \rangle = \sigma(N) = 0 \). Then, it remains 
\[
\langle \eta^S, \sigma \rangle \sigma = \left( \langle \mathbf{1}^S, \sigma \rangle - \frac{\lvert S \rvert}{\lvert N \rvert} \langle \mathbf{1}^N, \sigma \rangle \right) \sigma = \langle \mathbf{1}^S, \sigma \rangle \sigma = \sigma(S) \cdot \sigma.
\] 
\endproof

Because the size of the Gram matrices depends on the size of the collections and not on the number of players, we need to work with the smallest collections. Moreover, we see on Figure~\ref{fig: multi-proj} that the projection onto the core is not necessarily the projection associated with the largest independent collection. 

\begin{definition}
Let \( x \) be a preimputation outside of the core. We say that a collection \( \mathcal{Q} \) is a \emph{reaching collection} for \( x \) if it is independent and \( \pi_{A_\mathcal{Q}}(x) \in C(v) \). 
\end{definition}

We denote by \( \rho_\mathcal{Q}(x) \) the set of reaching subcollections of \( \mathcal{Q} \) for \( x \). If \( \mathcal{Q} \) is such that \( x \in X_\mathcal{Q}(v) \), we simply write \( \rho_\mathcal{Q}(x) \eqqcolon \rho(x) \). 

\begin{proposition}
Let \( x \in X(v) \setminus C(v) \). The projection of \( x \) onto the core coincide with the projection onto \( A_\mathcal{T} \), with \( \mathcal{T} \) being minimal in \( \rho(x) \).
\end{proposition}

\proof
For any face \( F \) of the core, we denote by \( \mathcal{E}_F \) the set of coalitions \( \mathcal{E}_F \coloneqq \{S \in \mathcal{N} \setminus \mathcal{E}(v) \mid F \subseteq A_S(v) \} \). Let \( x \in X_\mathcal{Q} \). The projection of \( x \) onto the core is necessarily achieved via a reaching collection. Let \( \mathcal{T}_0 \) be one of them. Then \( \pi_{A_{\mathcal{T}_0}}(x) \) belongs to a face of the core, denoted by \( F_0 \). If there does not exist another face \( F_1 \) containing \( F_0 \) such that \( \mathcal{E}_{F_1} \) belongs to \( \rho(x) \), then \( \mathcal{T}_0 \) is minimal in \( \rho(x) \). Indeed, if \( F_1 \supseteq F_0 \), then \( F_1 \) is included in less hyperplanes than \( F_0 \) and \( F_0 \) is included in all the hyperplanes in which \( F_1 \) lies. 

\medskip 

Then, we assume that there exists a face \( F_1 \supseteq F_0 \) such that \( \mathcal{E}_{F_1} \in \rho(x) \), and repeat the same process with \( \mathcal{E}_{F_1} \). Because the number of players is finite, the number of collection of coalitions is also finite, therefore there exists a face \( F_k \) containing \( F_0 \) such that \( \mathcal{E}_{F_k} \) is minimal in \( \rho(x) \). Because \( \pi_{A_{\mathcal{T}_0}}(x) \), which is the projection onto the core, belongs to \( F_0 \), it also belongs to \( F_k \). 
\endproof

\begin{algorithm}
\caption{Finding a minimal reaching collection} \label{algo: find-reach-coll}
\begin{algorithmic}[1]
\Require A balanced game \( (N, v) \), a preimputation \( x \)
\Ensure A minimal reaching collection
\State Find \( \mathcal{Q} \) such that \( x \in X_\mathcal{Q}(v) \)
\State Initialize \( \mathcal{K} \gets \{\{S\} \mid S \in \mathcal{Q} \} \)
\For{\( \mathcal{T} \in \mathcal{K} \)}
\If{\( \max_{S} e_S(\pi_{A_\mathcal{T}}(x)) \leq 0 \)}
\State \textbf{return} \( \mathcal{T} \)
\Else
\For{\( S \in \mathcal{Q} \)}
\If{\( \mathcal{T} \cup \{S\} \not \in \mathcal{K} \) and \Call{UpdateGramian}{$ \mathcal{T}, S $} \( \neq 0 \)}
\State Add \( \mathcal{T} \cup \{S\} \) to \( \mathcal{K} \)
\EndIf
\EndFor
\State Remove \( \mathcal{T} \) from \( \mathcal{K} \)
\EndIf
\EndFor
\end{algorithmic}
\end{algorithm}

We know that Algorithm \ref{algo: find-reach-coll} gives at least one output when the game is balanced, because the core has at least a nonempty face. 

\section{Application: measuring and correcting market failures}\label{sec: failure}

In the sequel, we propose a quantification of the \emph{failure} of a market, for a given payment vector, using the model of \emph{market games} defined by \textcite{shapley1969market}. A \emph{market} is a mathematical model, denoted by \( \mathfrak{m} = (N, G, A, U)\), where 
\begin{itemize}
\item \( N \) is a finite set of \emph{players}, or \emph{traders}, 
\item \( G \) is the nonnegative orthant of a finite-dimensional vector space, called the \emph{commodity space}, 
\item \( A = \{a^i \mid i \in N\} \subseteq G^N \) the elements of which are called \emph{initial endowments}, 
\item \( U = \{u^i \mid i \in N \} \) is an indexed set of continuous, concave function \( u^i: G \to \mathbb{R} \), called the \emph{utility functions}. Notice that there is no assumption concerning the monotonicity of the utility functions. 
\end{itemize}

The initial endowments can be interpreted as the belongings of each player at the initial state of the economy, and they will trade these commodities with other players to maximize their utility. Let \( S \) be a coalition. If the players in \( S \) form a market, the whole quantity of commodities is \( \sum_{i \in S} a^i \). Then, a set of endowments \( \{x^i \mid i \in S\} \subseteq G \) such that \( \sum_{i \in S} x^i = \sum_{i \in S} a^i \) is called a \emph{\( S \)-feasible allocation} of the market \( \mathfrak{m} \), and we denote their set by \( X^S \). The \emph{market game} generated by this market is a game \( (N, v) \) whose coalition function is defined, for all \( S \in \mathcal{N} \), by 
\[
v(S) = \max_{x \in X^S} \sum_{i \in S} u^i(x^i). 
\]
\textcite{shapley1969market} proved that market games are totally balanced, hence we can always project any preimputation onto it. The core of a market game is referred to as the core of the associated market. The core of a market is then the set of payment vectors that coalitions cannot achieve by forming a market on their own. 

\medskip

It is well-known that the competitive equilibria of an economy lie in its core. Consistently, \textcite{shapley1975competitive} proved that the set of competitive outcomes of a market is included in the core of its associated market game. Therefore, they satisfy very desirable properties, such as Pareto optimality and coalitional rationality, but the players of the market do not necessarily reach one of these points as a payment vector. Following the interpretation we adopted throughout this manuscript, a positive outcome for the interaction between players leading to the existence of the game must be a payment vector from the core, to avoid the departure of some players from the grand coalition. This leads to the following definition. 

\begin{definition}
Let \( \mathfrak{m} \) be a market, and \( (N, v) \) its associated game. Let \( x \) be a preimputation which is not included in the core. We define the \emph{failure} of the market \( \mathfrak{m} \), or of the game \( (N, v) \), at \( x \), denoted by \( \mu_v(x) \) or \( \mu_\mathfrak{m}(x) \), as the following quantity
\[
\mu_v(x) \coloneqq \min_{y \in \partial C(v)} \lVert x - y \rVert, 
\]
where \( \partial C(v) \) denotes the boundary of the core \( C(v) \). 
\end{definition}

The failure of the market at a specific preimputation \( x \) is the distance from \( x \) to the core. It represents the volume of the smallest reallocation of money, i.e., the Euclidean norm of a shortest side payment, needed to go from \( x \) to its closest element in the core. 

\medskip 

When players join a market, they put in common all their initial endowments, that they will share later. Players in \( S \), with the sum of their initial endowments, are able to redistribute them in a way that the sum of their utility after the redistribution is \( v(S) \). Because market games are totally balanced, the larger the grand coalition is, the better should be the payment of every coalition.

\medskip 

Then, when an allocation of commodities leads to a preimputation not included in the core, for at least one coalition \( S \) we have \( x(S) < v(S) \), which is absurd. Indeed, by definition there exists an allocation of the commodities of the players in the coalition \( S \) such that the sum of the utilities of these players is \( v(S) \). If there is a broader range of commodities to distribute, the payments of all the coalitions must be greater than the payment then can have by themselves. The failure at a specific preimputation \( x \) quantify how much \( x \) misses the core and fails to distribute correctly the payments. The players in \( S \) may want to leave this market, together with their initial endowments, which decreases the other players' payments because of the balancedness of the game. 

\medskip

The same interpretation holds for linear production games \cite{owen1975core} or flow games \cite{kalai1982totally}. When several players forming a coalition \( S \) join the game, they bring raw materials in the case of a linear production games, or new edges in the graph for a flow game. By themselves, they can achieve \( v(S) \), but by joining their effort with other players, for instance to have more diversity of raw materials, or a denser graph, with a greater capacity, we can expect that everyone is doing better than before joining, which can be an interpretation of being totally balanced. They cannot accept to have a payment lower than \( v(S) \), and the failure at a given preimputation quantifies the cumulative and nested inability of the market to give a proper payment to all the coalitions not satisfied by their payment. 

\begin{proposition}
Let \( \mathfrak{m} \) be a market, and \( (N, v) \) be its associated game. Let \( x \) be a preimputation not included in the core. Then
\[
\mu_v(x) = \min_{\mathcal{T} \in \rho(x)} \left\lVert \sum_{S \in \mathcal{T}} \gamma^\mathcal{T}_S(x) \eta^S \right\rVert = \min_{\mathcal{T} \in \rho(x)} \left\lVert \sum_{S \in \mathcal{T}} e_S(x) h^S \right\rVert. 
\]
\end{proposition}

\proof
Because the core is bounded, closed and convex, the projection of \( x \in X(v) \setminus C(v) \) lies on the frontier \( \partial C(v) \). By the Hilbert projection theorem, the projection is the point minimizing the distance, then 
\[
\mu_v(x) = \lVert x - \pi_{C(v)}(x) \rVert
\]
By the same argument, the projection onto the core coincides with the projection onto a face of the core that is the closest from the projected preimputation. Then, 
\[
\mu_v(x) = \min_{\mathcal{T} \in \rho(x)} \lVert x - \pi_{A_\mathcal{T}}(x) \rVert.
\]
Applying the formula of Theorem \ref{th: main-formula-proj} finishes the proof. 
\endproof

It may be possible to use the orthonormal basis of \( \mathcal{T} \) computed by Algorithm~\ref{algo: gramian} to let the norm and the sum commute. 

\medskip

Using the projectors to compute a market failure at a given preimputation also gives an explicit side payment to redistribute the global worth without changing it. Moreover, this side payment is the one of smallest norm. We can interpret it as the side payment with the smallest redistribution cost, and therefore the optimal reallocation from \( x \) to a core element. Let \( \mathcal{T} \) be the reaching collection defining the core projection. Then the optimal reallocation 
\[
\sigma_{\text{opt}} = \sum_{S \in \mathcal{T}} \gamma^\mathcal{T}_S(x) \eta^S
\]
works as follows: for each coalition \( S \in \mathcal{T} \), we collect from each player the worth \( \frac{\lvert S \rvert}{\lvert N \rvert} \gamma^\mathcal{T}_S(x)\), and give \( \gamma^\mathcal{T}_S(x) \) to each player in \( S \).  

\medskip 

The preimputations belonging to the frontier of the core have a failure of \( 0 \), and we can even extend this function to an element \( x \) of the interior of the core by \( \mu_v(x) \coloneqq - \min_{y \in \partial C(v)} \lVert x - y \rVert \), and merge these two functions, for all \( x \in X(v) \), by 
\[
\mu_v(x) = (-1)^{\mathbf{1}_{C(v)}(x)} \min_{y \in \partial C(v)} \lVert x - y \rVert, 
\] 
with \( \mathbf{1}_{C(v)} \) being the indicator function of the core. These quantities are linked to solutions concepts developed by the same authors, before defining market games. To read more about it, see \textcite[][Chapter 7]{peleg2007introduction}. We denote \( \mathcal{N}_* \coloneqq \mathcal{N} \setminus \{N\} \). 

\begin{definition}[\textcite{shapley1966quasi}] \leavevmode \newline
Let \( \varepsilon \) be a real number. The \( \varepsilon \)-core of the game \( (N, v) \), denoted by \( C_\varepsilon(v) \), is defined by 
\[
C_\varepsilon(v) \coloneqq \{x \in X(v) \mid e_S(x) \leq \varepsilon, \hspace{1pt} \forall S \in \mathcal{N}_* \}. 
\]
\end{definition}

Notice that \( C_0(v) = C(v) \). They also defined the smallest of these sets, called the \emph{least-core} and denoted it by \( LC(v) \), as the intersection of all nonempty \( \varepsilon \)-cores of \( (N, v) \). It is possible to define the least-core as an \( \varepsilon_0 \)-core, with
\[ 
\varepsilon_0 = \min_{x \in X(v)} \max_{S \in \mathcal{N}_*} \, e_S(x).
\] 
Then, we have that 
\[
LC(v) = \argmin_{x \in X(v)} \max_{S \in \mathcal{N}_*} \, e_S(x). 
\] 

\begin{definition}
Let \( \mathfrak{m} \) be a market, and \( (N, v) \) its associated game. We define the \emph{Chebyshev core}, denoted by \( C_\mu (v) \), by 
\[
C_\mu (v) = \argmin_{x \in X(v)} \mu_v(x). 
\]
\end{definition}

The projection of an element from the interior of a polytope onto its frontier is the projection onto a facet of the polytope. Indeed, any face of the polytope is the intersection of some facets, and by the triangle inequality, the distance from an intersection of facets is longer than the distance from a facet. Then, the distance between an element \( x \in C(v) \) and \( \partial C(v) \) is given by 
\[ 
\min_{y \in \partial C(v)} \lVert x - y \rVert = \min_{S \in \mathcal{N}_*} \, \lVert x - \pi_{A_S}(x) \rVert.
\] 
Using the formula of Proposition \ref{prop: first-proj} leads to
\[
\min_{y \in \partial C(v)} \lVert x - y \rVert = \min_{S \in \mathcal{N}_*} \, \left\lVert \frac{e_S(x)}{\lVert \eta^S \rVert^2} \eta^S \right\rVert = \min_{S \in \mathcal{N}_*} \, \frac{\lvert e_S(x) \rvert}{\lVert \eta^S \rVert}. 
\]
Because \( x \) belongs to the core, its excess is nonpositive, then 
\[
\min_{y \in \partial C(v)} \lVert x - y \rVert = \max_{S \in \mathcal{N}_*} \frac{e_S(x)}{\lVert \eta^S \rVert} = \max_{S \in \mathcal{N}_*} \gamma_S(x),
\]
and we can write 
\[
C_\mu (v) = \argmin_{x \in X(v)} \max_{S \in \mathcal{N}_*} \, \gamma_S(x), 
\]
which resembles \( LC(v) \). Similarly to the least-core, the Chebyshev core is never empty. 

\medskip 

The difference between the least-core and the Chebyshev core is the rescaling factor \( \lVert \eta^S \rVert^{-1} \). It can be interpreted as a coefficient scaling the excess at the individual level, because \( \lVert \eta^S \rVert \) only depends on the cardinality of \( S \). If the excess of a \( 5 \)-player coalition \( S \) is \( 4 \), and the excess of a \( 2 \)-player coalition \( T \) is \( 2 \), despite the fact that the excess for \( S \) is bigger than for \( T \), the players in \( T \) are more aggrieved than players in \( S \). The average excess per player in \( T \) is \( 1 \), while the average excess per player in \( S \) is \( \frac{4}{5} < 1 \). 

\medskip 

Finally, the Chebyshev core of a market, or, equivalently, of a game, can be efficiently computed using a unique linear program \cite{boyd2004convex}. The number of constraints for this program can be reduced using the algorithms by the author, Grabisch and Sudh{\"o}lter~\cite{laplace2023minimal}. 

\section{Concluding remarks}

In this paper, we have presented formulae and algorithms to project a preimputation onto the core. The main difficulty in this task is the exponential size of the H-description of the core, as well as the impossibility to know beforehand which inequalities are redundant in the description of the core. 

\medskip 

The projector onto the core is built from two parts. First, an exact formula giving the projection onto specific affine subspaces of preimputations (Section~\ref{sec: projection-arbitrary-intersection}) is provided. Second, an algorithm that identifies in which face the projection lies, and which set of coalitions corresponds to this face (Section~\ref{sec: algo-proj}) is designed. 

\medskip 

Because we know which core allocation is the closest from the current state of the economy, we are able to compute the distance from this state to the core, and then to compute a \emph{failure} map (Section~\ref{sec: failure}). Because a positive failure represents how inefficient is a given preimputation to redistribute the aggregated endowments of the players, a negative failure represents a net increase in the overall well-being. Following this idea, we have defined a new solution concept, the \emph{Chebyshev core}, which maximizes this social welfare and only consists of socially-optimal outcomes. 

\printbibliography

\end{document}